\documentclass[aos,nameyear,seceqn,rotating,dvips]{arximspdf}
\usepackage{dcolumn}
\usepackage{graphicx}

\doi{10.1214/10-AOS832}
\volume{39}
\issue{1}
\pubyear{2011}
\firstpage{174}
\lastpage{200}

\makeatletter
\newcolumntype{d}[1]{D{.}{.}{#1}}
\newcommand{\es}{S}
\newcommand{\betac}{\beta_{c}}
\newcommand{\betau}{\beta_{u}}
\newcommand{\betaczero}{\beta_{c,0}}
\newcommand{\betauS}{\beta_{u,\es}}
\newcommand{\dc}{{d_{c}}}
\newcommand{\du}{{d_{u}}}
\newcommand{\duS}{{d_{u,\es}}}

\newtheorem{theorem}{Theorem}%[section]
\newtheorem{lemma}{Lemma}[section]
\newproclaim{remark}{Remark}%[section]     %
\newcommand{\diag}{\mathrm{diag}}

\newcommand{\T}{\mathrm{T}}

\newcommand{\bfB}{\mathbf{B}}
  \newcommand{\bfD}{\mathbf{D}}

  \newcommand{\bfG}{\mathbf{G}}

\newcommand{\bft}{\mathbf{t}}

\newcommand{\bfx}{\mathbf{x}}  
  
\newcommand{\bfz}{\mathbf{z}}  
\makeatother

\begin{document}
\begin{frontmatter}

\title{Focused information criterion and model averaging for generalized additive partial linear models}
\runtitle{Generalized additive partially linear models}

\begin{aug}
\author[A]{\fnms{Xinyu} \snm{Zhang}\thanksref{t1}\ead[label=e1]{xinyu@amss.ac.cn}}
\and
\author[B]{\fnms{Hua} \snm{Liang}\thanksref{t2}\ead[label=e3]{hliang@bst.rochester.edu}\corref{}}

\thankstext{t1}{Supported in part by the National Natural Science Foundation of China Grants 70625004 and 70933003.}
\thankstext{t2}{Supported in part by NSF Grant DMS-08-06097.}
\runauthor{X. Zhang and H. Liang}
\affiliation{Chinese Academy of Sciences and University of Rochester}
\address[A]{Institute of Systems Science\\
Academy of Mathematics and System Science\\
Chinese Academy of Sciences\\
Beijing, 100190\\
China\\
\printead{e1}} %adresu isvedimo komanda gale!
\address[B]{Department of Biostatistics\\
\quad and Computational Biology\\
University of Rochester\\
Rochester, New York 14642\\
USA\\
\printead{e3}}
\end{aug}

% HISTORY:
\received{\smonth{2} \syear{2010}}
\revised{\smonth{5} \syear{2010}}

% ABSTRACT
\begin{abstract}
We study model selection and model averaging in generalized additive
partial linear models (GAPLMs). Polynomial spline  is used to
approximate nonparametric functions. The corresponding estimators of
the linear parameters are shown to be asymptotically normal. We then
develop a focused information criterion (FIC) and a frequentist model
average (FMA) estimator on the basis of the quasi-likelihood principle
and examine theoretical  properties of the FIC and FMA. The major
advantages of the proposed procedures over the existing ones are their
computational expediency and theoretical reliability. Simulation
experiments have provided evidence of the superiority of the proposed
procedures. The approach is further applied to a real-world data
example.
\end{abstract}



\end{frontmatter}

%s1 ###
\section{Introduction}\label{sec:intr}

Generalized additive models, which are a generalization of the generalized models and involve a summand of
one-dimensional nonparametric functions instead of a summand of linear components, have been widely used to
explore the complicated relationships between a response to treatment
and predictors of interest [Hastie and Tibshirani (\citeyear{ht1990})].
Various attempts are still being made to balance the interpretation of
generalized linear models and the flexibility of generalized additive
models such as generalized additive partial linear models (GAPLMs), in
which some of the additive component functions are linear, while the
remaining ones are modeled nonparametrically [H\"ardle et al.
(\citeyear{hetal2004a,hetal2004b})].
A special case of a GAPLM with a single nonparametric component, the
generalized partial linear model (GPLM), has been well studied in the
literature; see, for example, Severini and Staniswalis
(\citeyear{ss1994}), Lin and Carroll (\citeyear{lc2001}), Hunsberger
(\citeyear{h1994}), Hunsberger et al. (\citeyear{heta2002}) and Liang
(\citeyear{l2008}). The profile quasi-likelihood procedure has
generally been used, that is, the estimation of GPLM is made
computationally feasible by the idea that estimates of the parameters
can be found for a known nonparametric function, and an estimate of the
nonparametric function can be found for the estimated parameters.
Severini and Staniswalis (\citeyear{ss1994}) showed that the resulting
estimators of the parameter are asymptotically normal and that
estimators of the nonparametric functions are consistent in supremum
norm. The computational algorithm involves searching for maxima of
global and local likelihoods simultaneously. It is worthwhile to point
out that studying GPLM is easier than studying GAPLMs, partly because
there is only one nonparametric term in GPLM. Correspondingly,
implementation of the estimation for GPLM is simpler than for GAPLMs.
Nevertheless, the GAPLMs are more flexible and useful than GPLM because
the former allow several nonparametric terms for some covariates and
parametric terms for others, and thus it is possible to explore more
complex relationships between the response variables and covariates.
For example, Shiboski (\citeyear{s1998}) used a GAPLM to study AIDS
clinical trial data and M\"uller and R\"onz (\citeyear{mr2000}) used a
GAPLM to carry out credit scoring. However, few theoretical results are
available for GAPLMs, due to their general flexibility. In this
article, we shall study estimation of GAPLMs using polynomial spline,
establish asymptotic normality for the estimators of the linear
parameters and develop a focused information criterion (FIC) for model
selection and a frequentist model averaging (FMA) procedure in
construction of the confidence intervals for the focus parameters with
improved coverage probability.

We know that traditional model selection methods such as the Akaike
information criterion [AIC, Akaike (\citeyear{a1973})] and the Bayesian
information criterion [BIC, Schwarz (\citeyear{s1978})] aim to select a
model with good overall properties, but the selected model is not
necessarily good for estimating a specific parameter under
consideration, which may be a function of the model parameters; see an
inspiring example in Section 4.4 of Claeskens and Hjort
(\citeyear{ch2003}). Exploring the data set from the Wisconsin
epidemiologic study of diabetic retinopathy, Claeskens, Croux and van
Kerckhoven (\citeyear{ccvk2006}) also noted that different models are
suitable for different patient groups. This occurrence has been
confirmed by Hand and Vinciotti (\citeyear{hv2003}) and Hansen
(\citeyear{h2005}). Motivated by this concern, Claeskens and Hjort
(\citeyear{ch2003}) proposed a new model selection criterion, FIC,
which is an unbiased estimate of the limiting risk for the limit
distribution of an estimator of the focus parameter, and systematically
developed a general asymptotic theory for the proposed criterion.  More
recently, FIC has been studied in several models. Hjort and Claeskens
(\citeyear{hc2006}) developed the FIC for the Cox hazard regression
model and applied it to a study of skin cancer; Claeskens, Croux and
van Kerckhoven (\citeyear{ccvk2007}) introduced the FIC for
autoregressive models and used it to predict the net number of new
personal life insurance policies for a large insurance company.

The existing model selection methods may arrive at a model which is
thought to be able to capture the main information of the data, and to
be decided in advance in data analysis. Such an approach may lead to
the ignoring of  uncertainty introduced by model selection. Thus, the
reported confidence intervals are too narrow or shift away from the
correct location, and the corresponding coverage probabilities of the
resulting confidence intervals can substantially deviate from the
nominal level [Danilov and Magnus (\citeyear{dm2004}) and Shen, Huang
and Ye (\citeyear{shy2004})]. Model averaging, as an alternative to
model selection, not only provides a kind of insurance against
selecting a very poor model, but can also avoid model selection
instability [Yang (\citeyear{y2001}) and Leung and Barron
(\citeyear{lb2006})] by weighting/smoothing estimators across several
models, instead of relying entirely on a single model selected by some
model selection criterion. As a consequence, analysis of the
distribution of model averaging estimators can improve coverage
probabilities. This strategy has been adopted and studied in the
literature, for example, Draper (\citeyear{d1995}), Buckland, Burnham
and Augustin (\citeyear{bba1997}), Burnham and Anderson
(\citeyear{ba2002}), Danilov and Magnus (\citeyear{dm2004}) and Leeb
and P\"{o}stcher (\citeyear{lp2006}). A seminal work, Hjort and
Claeskens (\citeyear{hc2003}), developed asymptotic distribution
theories for estimation and inference after model selection and model
averaging across parametric models. See Claeskens and Hjort
(\citeyear{ch2008}) for a comprehensive survey on FIC and model
averaging.

FIC and FMA have been well studied for parametric models. However, few
efforts have been made to study FIC and FMA for semiparametric models.
To the best of our knowledge, only Claeskens and Carroll
(\citeyear{cc2007}) studied FMA in semiparametric partial linear models
with a univariate nonparametric component. The existing results are
hard to  extend directly to GAPLMs, for the following reasons: (i)
there exist nonparametric components in GAPLMs, so the ordinary
likelihood method cannot be directly used in estimation for GAPLMs;
(ii) unlike the semiparametric partial linear models in Claeskens and
Carroll (\citeyear{cc2007}), GAPLMs allow for multivariate covariate
consideration in nonparametric components and also allow for the mean
of the response variable to be connected to the covariates by a link
function, which means that the binary/count response variable can be
considered in the model. Thus, to develop FIC and FMA procedures for
GAPLMs and to establish asymptotic properties for these procedures are
by no means straightforward to achieve. Aiming at these two goals, we
first need to appropriately estimate the coefficients of the parametric
components (hereafter, we call these coefficients ``linear
parameters'').

There are two commonly used estimation approaches for GAPLMs: the first
is local scoring backfitting, proposed by Buja, Hastie and Tibshirani
(\citeyear{bht1989}); the second is an application of the marginal
integration approach on the nonparametric component [Linton and
Nielsen (\citeyear{ln1995})]. However, theoretical properties of the
former are not well understood since it is only defined implicitly as
the limit of a complicated iterative algorithm, while the latter
suffers from the \textit{curse of dimensionality} [H\"{a}rdle et al.
(\citeyear{hetal2004a})], which may lead to an increase in the
computational burden and which also conflicts with the purpose of using
a GAPLM, that is, dimension reduction. Therefore, in this article, we
apply polynomial spline to approximate nonparametric functions in
GAPLMs. After the spline basis is chosen, the nonparametric components
are replaced by a linear combination of spline basis, then the
coefficients can be estimated by an efficient one-step maximizing
procedure. Since the polynomial-spline-based method solves much smaller
systems of equations than kernel-based methods that solve larger
systems (which may lead to identifiability problems), our
polynomial-spline-based procedures can substantially reduce the
computational burden. See a similar discussion about this computational
issue in Yu, Park and Mammen (\citeyear{ypm2008}), in the generalized
additive models context.

The use of polynomial spline  in generalized nonparametric models can
be traced back to Stone (\citeyear{s1986}), where the rate of
convergence of the polynomial spline estimates for the generalized
additive model were first obtained. Stone (\citeyear{s1994}) and Huang
(\citeyear{h1998}) investigated the polynomial spline estimation for
the generalized functional ANOVA model. In a widely discussed paper,
Stone et al. (\citeyear{setal1997}) presented a completely theoretical
setting of polynomial spline approximation, with applications to a wide
array of statistical problems, ranging from least-squares regression,
density and conditional density estimation, and generalized regression
such as logistic and Poisson regression, to polychotomous regression
and hazard regression. Recently, Xue and Yang (\citeyear{xy2006})
studied estimation in the additive coefficient model with continuous
response using polynomial spline
 to approximate the coefficient functions. Sun, Kopciuk and Lu (\citeyear{skl2008}) used
polynomial spline in partially linear single-index proportional hazards
regression models.  Fan, Feng and Song (\citeyear{ffs2009}) applied
polynomial spline to develop nonparametric independence screening in
sparse ultra-high-dimensional additive models. Few attempts have been
made to study polynomial spline for GAPLMs, due to the extreme
technical difficulties involved.

The remainder of this article is organized as follows. Section \ref{SEC:framework and estimation} sets
out the model framework and provides the polynomial spline estimation
and asymptotic normality of estimators. Section \ref{SEC:model selection} introduces the FIC
and FMA procedures and constructs confidence intervals for the focus
parameters on a basis of FMA estimators. A simulation study and
real-world data analysis are presented in Sections \ref{s4} and \ref{s5},
respectively. Regularity conditions and technical proofs are presented
in the \hyperref[appe]{Appendix}.

%s2 ###
\section{Model framework and estimation}\label{SEC:framework and estimation}

We consider a GAPLM where the response $Y$ is related to covariates
$\mathbf{X}=( X_{1},\ldots,X_{p})^{\T}\in R^{{p}}$ and $\mathbf{Z}=(
Z_{1},\ldots,\break Z_{d})^{\T}\in R^{d}$. Let the unknown mean response $
\mathbf{u} ( \mathbf{x},\mathbf{z}) =E( Y|\mathbf{X}=\mathbf{x},
\mathbf{Z}=\mathbf{z}) $ and the conditional variance function be
defined by a known positive function $V$, $
\operatorname{var}(Y|\mathbf{X}=\mathbf{x},\mathbf{Z}=\mathbf{z})
=V\{\mathbf{u}(\mathbf{x},\mathbf{z})\} $.
In this article, the mean function $\mathbf{u}$ is defined via a known
link function $g$ by an additive linear function
%e2.1 ###
\begin{equation}
g\{ \mathbf{u} ( \mathbf{x},\mathbf{z}) \} =\sum_{\alpha=1}^{{p}}\eta
_{\alpha}( x_{\alpha})+\mathbf{z}^{\T}\beta , \label{model:GAPLM}
\end{equation}
where $x_\alpha$ is the $\alpha$th element of $x$, $\bolds{\beta }$ is
a $d$-dimensional regression parameter and the $\eta _{\alpha}$'s are
unknown smooth functions. To ensure identifiability, we assume that
$E\{\eta _{\alpha }(X_{\alpha })\}=0$ for $1\leq \alpha \leq {p}$.

Let $\beta=(\betac^{\T}, \betau^{\T} )^{\T}$ be a vector with
$d=\dc+\du$ components, where $\betac$ consists of the first $\dc$
parameters of $\beta$ (which we certainly wish to be in the selected
model) and $\betau$ consists of the remaining $\du$ parameters (for
which we are unsure whether or not they should be included in the
selected model). In what follows, we call the elements of $\bfz$
corresponding to $\betac$ and $\betau$ the \textit{certain} and
\textit{exploratory} variables, respectively. As in the literature on
FIC, we consider a local misspecification framework where the true
value of the parameter vector $\beta$ is
$\beta_0=(\betaczero^{\T},\delta^{\T}/\sqrt{n})^{\T}$, with $\delta$
being a $\du\times 1$ vector; that is, the true model is away from the
deduced model with a distance $O(1/{\sqrt n})$. This framework
indicates that squared model biases and estimator variances are both of
size $O(1/n)$, the most possible large-sample approximations. Some
arguments related to this framework appear in Hjort and Claeskens
(\citeyear{hc2003,hc2006}).

Denote by $\beta_{\es}=(\betac^{\T}, \betauS^{\T})^{\T}$ the parameter
vector in the ${\es}$th submodel, in the same sense as $\beta$, with
$\betauS$ being a $\duS$-subvector of $\beta_u$. Let $\pi_{\es}$ be the
projection matrix of size $\duS\times d_{u}$ mapping $\beta_{u}$ to
$\betauS$. With $\du$ exploratory covariates, our setup allows
$2^{\du}$ extended models to choose among. However, it is not necessary
to deal with all $2^{\du}$ possible models and one is free to consider
only a few relevant submodels (unnecessarily nested or ordered) to be
used in the model selection or averaging. A special example is the
James--Stein-type estimator studied by Kim and White
(\citeyear{kw2001}), which is a weighted summand of the estimators
based on the reduced model ($\duS=0$) and the full model ($\duS=\du$).
So, the covariates in the ${\es}$th submodel are $\mathbf{X}$ and
$\Pi_{\es}\mathbf{Z}$, where $\Pi_{\es}=\diag(I_{\dc}, \pi_{\es}).$ To
save space, we generally ignore the dimensions of zero vectors/matrices
and identity matrices, simply denoting them by $0$  and $I,$
respectively. If necessary, we will write their dimensions explicitly.
In the remainder of this section, we shall investigate polynomial
spline estimation for $(\betaczero^{\T}, 0)$ based on the ${\es}$th
submodel and establish a theoretical property for the resulting
estimators.

Let $\eta _{0}=\sum_{\alpha=1}^{{p}}\eta _{0,\alpha}( x_{\alpha})$ be
the true additive function and the covariate $X_{\alpha }$ be
distributed on a compact interval $[a_{\alpha },b_{\alpha }]$. Without
loss of generality, we take all intervals $[a_{\alpha },b_{\alpha
}]=[0,1]$ for $\alpha =1,\ldots,{p}$. Noting (\ref{eq:etainfty}) in
Appendix~\ref{s62}, under some smoothness assumptions in Appendix \ref{s61}, $\eta
_{0}$ can be well approximated by spline functions. Let
$\mathcal{S}_{n}$ be the space of polynomial splines on $[0,1]$ of
degree \mbox{$\varrho \geq 1$}. We introduce a knot sequence with $J$ interior
knots, $ k_{-\varrho
}=\cdots=k_{-1}=k_{0}=0<k_{1}<\cdots<k_{J}<1=k_{J+1}=\cdots=k_{J+\varrho
+1}, $ where $J\equiv J_{n}$ increases when sample size $n$ increases
and the precise order is given in condition (C6). Then,
$\mathcal{S}_{n}$ consists of functions $\varsigma$ satisfying the following:
\begin{enumerate}[(ii)]
\item[(i)] $\varsigma$ is a polynomial of degree $\varrho $ on each of the subintervals $[ k_{j},k_{j+1}) $, $j=0,\ldots,J_{n}-1$, and the last subinterval is $[ k_{J_{n}},1] $;
\item[(ii)] for $\varrho \geq 2$, $\varsigma$ is ($\varrho -1$)-times continuously differentiable on $[0,1]$.
\end{enumerate}
For simplicity of proof, equally spaced knots are used. Let $h=1/(J_{n}+1) $ be the distance between two consecutive knots.

Let $( Y_{i},\mathbf{X}_{i},\mathbf{Z}_{i}) $, $i=1,\ldots,n$, be
independent copies of $( Y,\mathbf{X},\mathbf{Z}) $. In the ${\es}$th
submodel, we consider the additive spline estimates of $\eta
_{0} $ based on the independent random sample $( Y_{i},\mathbf{X}_{i},%
\Pi_{\es}\mathbf{Z}_{i}) $, $i=1,\ldots,n$. Let $\mathcal{G}_{n}$ be the
collection of functions $\eta $ with the additive form {$\eta ( \mathbf{x
}) =\sum_{\alpha=1}^{{p}}\eta _{\alpha}(
x_{\alpha})$}, where each component function $\eta _{\alpha }\in \mathcal{S}
_{n}$.

We would like to find a function $\eta \in \mathcal{G}_{n}$ and a value of $
\bolds{\beta}_{\es}$ that maximize the  quasi-likelihood function
%e2.2 ###
\begin{equation}\label{DEF:quasilikelihood}
L(\eta,\bolds{\beta}_{\es})
=
\frac{1}{n}\sum_{i=1}^{n}Q[g^{-1}\{\eta(\mathbf{X}_{i})+(\Pi_{\es}\mathbf{Z}_{i})^{\T}\bolds{\beta}_{\es}\} ,Y_{i}],
\qquad
\eta\in\mathcal{G}_{n}\mathbf,
\end{equation}
where $Q( m,y)$ is the quasi-likelihood function satisfying
$\frac{\partial Q(m,y)}{\partial m}=\frac{y-m}{V(m)}$.

For the $\alpha$th covariate $x_{\alpha }$, let $b_{j,\alpha}(x_{\alpha})$ be the B-spline basis function of degree~$\varrho$. For
any $\eta\in\mathcal{G}_{n}$, one can write
$\eta(\mathbf{x})=\bolds{\gamma }^{\T}\mathbf{b}(\mathbf{x})$, where
$\mathbf{b}(\mathbf{x}) =\{b_{j,\alpha}(x_{\alpha}), j=-\varrho,\ldots,J_{n}, {\alpha =1,\ldots,{p}}\}^{\T}$ are the spline basis functions
and $\bolds{\gamma}=\{\gamma_{j,\alpha },j=-\varrho,\ldots,J_{n},\alpha =1,\ldots,{p}\}^{\T}$ is the spline coefficient vector. Thus, the
maximization problem in (\ref{DEF:quasilikelihood}) is equivalent to
finding values of $\bolds{\beta}_{\es}^{\ast}$ and
$\bolds{\gamma}^{\ast}$ that maximize
%e2.3 ###
\begin{equation}\label{DEF:quasilikelihood1}
\frac{1}{n}\sum_{i=1}^{n}Q[g^{-1}\{\bolds{\gamma}^{{\ast}\mathrm{T}}\mathbf{b}(\mathbf{X}_{i})+(\Pi_{\es}\mathbf{Z}_{i})^{\T}\bolds{\beta }_{\es}^{\ast}\},Y_{i}].
\end{equation}
We denote the maximizers as $\widehat{\bolds{\beta}}{}^{\ast}_{\es}$ and $\widehat{\bolds{\gamma}}{}^{\ast}_{\es}=
\{\widehat{\gamma}{}^{\ast}_{\es,j,\alpha },j=-\varrho ,\ldots,J_{n},\alpha =1,\ldots,\break{p}\}^{\T}$. The spline estimator of $\eta _{0}$ is then
$\widehat{\eta}{}^{\ast}_{\es}=\widehat{\bolds{\gamma}}{}_{\es}^{{\ast}\T}\mathbf{b}(\mathbf{x})$ and the centered spline estimators of each
component function are
\[
\widehat{\eta}{}^{\ast}_{\es,\alpha}(x_{\alpha})=\sum_{j=-\varrho}^{J_{n}}\widehat{\gamma}{}^{\ast}_{\es,j,\alpha}b_{j,\alpha}( x_{\alpha })-\frac{1}{n}\sum_{i=1}^{n}\sum_{j=-\varrho}^{J_{n}}\widehat{\gamma}{}^{\ast}_{\es,j,\alpha}b_{j,\alpha}(X_{i\alpha}),\qquad\alpha=1,\ldots,{p}.
\]
The above estimation approach can be easily implemented with commonly
used statistical software since the resulting model is a generalized
linear model.

For any measurable functions $\varphi_{1}$, $\varphi_{2}$ on $[0,1]^{{p}}$, define the empirical inner product and the corresponding norm
as
\[
\langle\varphi_{1},\varphi_{2}\rangle_{n}=n^{-1}\sum_{i=1}^{n}\{\varphi_{1}(\mathbf{X}_{i})\varphi_{2}(\mathbf{X}_{i})\},
\qquad
\|\varphi\|_{n}^{2}=n^{-1}\sum_{i=1}^{n}\varphi^{2}(\mathbf{X}_{i}).
\]
If $\varphi_{1}$ and $\varphi_{2}$ are $L^{2}$-integrable, define the
theoretical inner product and the corresponding norm as $\langle\varphi_{1},\varphi_{2}\rangle=E\{\varphi_{1}(\mathbf{X})\varphi_{2}(\mathbf{X})\}$,
$\|\varphi\|_{2}^{2}=E\varphi^{2}(\mathbf{X})$, respectively.
Let $\|\varphi\|_{n\alpha}^{2}$ and $\|\varphi\|_{2\alpha}^{2}$ be the empirical and theoretical norms, respectively, of a function $\varphi$ on $[0,1]$, that is,
\[
\|\varphi\|_{n\alpha}^{2}=n^{-1}\sum_{i=1}^{n}\varphi^{2}(X_{i\alpha }),\qquad\|\varphi\| _{2\alpha}^{2}=E\varphi^{2}(X_{\alpha})=\int_{0}^{1}\varphi^{2}( x_{\alpha})f_{\alpha}(x_{\alpha})\,dx_{\alpha },
\]
where $f_{\alpha}(x_{\alpha})$ is the density function of $X_{\alpha}$.

Define the  centered version spline basis for any $\alpha =1,\ldots,{p}$
and $j=-\varrho+1,\ldots,J_{n}$, $b_{j,\alpha}^{*}(x_{\alpha})=b_{j,\alpha}(x_{\alpha})-{\|b_{j,\alpha}\|_{2\alpha}}/{\|b_{j-1,\alpha}\|_{2\alpha}}b_{j-1,\alpha}(x_{\alpha})$, with the
standardized version given by
%e2.4 ###
\begin{equation}\label{DEF:BJalpha}
B_{j,\alpha}(x_{\alpha})=\frac{b^{\ast}_{j,\alpha}(x_{\alpha})}{\|b^{\ast}_{j,\alpha}\|_{2\alpha}}.
\end{equation}
Note that to find $(\gamma^{\ast},\beta_{\es}^{\ast})$ that maximizes (\ref{DEF:quasilikelihood1})
is mathematically equivalent to finding $(\gamma,\beta_{\es})$ that maximizes
%e2.5 ###
\begin{equation}\label{DEF:quasilikelihood2}
\ell(\bolds{\gamma},\bolds{\beta}_{\es})=\frac{1}{n}\sum_{i=1}^{n}Q[ g^{-1}\{\bolds{\gamma}^{\T}\mathbf{B}(\mathbf{X}_{i})+(\Pi_{\es}\mathbf{Z}_{i})^{\T}\bolds{\beta}_{\es}\},Y_{i}],
\end{equation}
where $\mathbf{B}(\mathbf{x})=\{B_{j,\alpha}(x_{\alpha}),j=-\varrho+1,\ldots,J_{n},{\alpha =1,\ldots,{p}}\}^{\T}$.
Similarly to $\widehat{\bolds{\beta}}{}_{\es}^{\ast}$, $\widehat{\bolds{\gamma}}{}_{\es}^{\ast}$, $\widehat{\eta}{}_{\es}^{\ast}$ and
$\widehat{\eta}{}_{\es,\alpha}^{\ast}$, we can define
$\widehat{\bolds{\beta}}_{\es}$, $\widehat{\bolds{\gamma}}_{\es}$,
$\widehat{\eta}_{\es}$ and the centered spline estimators of each
component function $\widehat{\eta}_{\es,\alpha}(x_\alpha)$.
In practice, the basis $\{b_{j,\alpha}(x_{\alpha}),j=-\varrho,\ldots,J_{n}, \alpha =1,\ldots,p\}^{\T}$
is used for data analytic implementation and the mathematically equivalent expression (\ref{DEF:BJalpha}) is convenient for asymptotic derivation.

Let $\rho_{l}(m)=\{\frac{dg^{-1}(m)}{dm}\}^l/V\{g^{-1}(m)\}$,
$l=1,2$. {Write $\mathbf{T}=(\mathbf{X}^{\T},\mathbf{Z}^{\T})^{\T}$, $m_{0}(\mathbf{T})=\eta_{0}(\mathbf{X})+\mathbf{Z}^{\T}\bolds{\beta}_{0}$
and $\varepsilon=Y-g^{-1}\{m_{0}(\mathbf{T})\}$. $\mathbf{T}_i$, $m_{0}(\mathbf{T}_i)$ and
$\varepsilon_i$ are defined in the same way after replacing
$\mathbf{X}$, $\mathbf{Z}$ and $\mathbf{T}$ by $\mathbf{X}_i$,
$\mathbf{Z}_i$ and $\mathbf{T}_i$, respectively.
Write
\begin{eqnarray*}
\bolds{\Gamma}(\mathbf{x})
&=&
\frac{E[\mathbf{Z}\rho_{1}\{m_{0}(\mathbf{T})\}|\mathbf{X}=\mathbf{x}]}{E[\rho_{1}\{m_{0}(\mathbf{T})\}|\mathbf{X}=\mathbf{x}]},\qquad\psi(\mathbf{T})=\mathbf{Z}-\bolds{\Gamma}(\mathbf{X}),
\\
\mathbf{G}_n
&=&
\frac{1}{\sqrt{n}}\sum_{i=1}^{n}\varepsilon_i\rho_1\{m_0(\mathbf{T}_{i})\}\bolds{\psi}(\mathbf{T}_{i}),\qquad\mathbf{D}=E[\rho_{1}\{m_0(\mathbf{T})\}\psi(\mathbf{T})\{\psi(\mathbf{T})\}^{\T}]
\end{eqnarray*}
 and
$\bolds{\Sigma}=E[\rho_{1}^{2}\{m_0(\mathbf{T})\}\varepsilon^{2}\psi(\mathbf{T})\{\psi(\mathbf{T})\}^{\T}]$.

The following theorem shows that the estimators $\widehat\beta_{\es}$
on the basis of the ${\es}$th submodel are asymptotically normal.

\begin{theorem}\label{th1}
Under the local misspecification framework and conditions \textup{(C1)--(C11)} in the \hyperref[appe]{Appendix},
\begin{eqnarray*}
&&\sqrt{n}\{\widehat\beta_{\es}-( \betaczero^\T, 0)^\T\}
\\
&&\qquad=
-(\Pi_{\es}\mathbf{D}\Pi_{\es}^{\T})^{-1}\Pi_{\es}\mathbf{G}_n+(\Pi_{\es}\mathbf{D}\Pi_{\es}^{\T})^{-1}\Pi_{\es}\mathbf{D}
\pmatrix{
0\cr
\delta}+o_p(1)
\\
&&\qquad\stackrel{d}\longrightarrow
-(\Pi_{\es}\mathbf{D}\Pi_{\es}^{\T})^{-1}\Pi_{\es}\mathbf{G}+(\Pi_{\es}\mathbf{D}\Pi_{\es}^{\T})^{-1}\Pi_{\es}\mathbf{D}
\pmatrix{
0\cr
\delta}
\end{eqnarray*}
with $\mathbf{G}_n\stackrel{d}\longrightarrow\mathbf{G}\sim N(0,\bolds{\Sigma})$,  where ``$\stackrel{d}\longrightarrow$'' denotes
convergence in distribution.
\end{theorem}

\begin{remark}\label{re:1}
If the link function $g$ is identical and there is only one
nonparametric component (i.e., $p=1$), then the result of Theorem \ref{th1} will simplify to those of Theorems 3.1--3.4  of Claeskens and
Carroll (\citeyear{cc2007}) under the corresponding submodels.
\end{remark}

\begin{remark}\label{re:2}
Assume that $\du=0$. Theorem \ref{th1} indicates that the
polynomial-spline-based estimators of the linear parameters are
asymptotically normal. This is the first explicitly theoretical result
on  asymptotic normality for estimation of the linear parameters in
GAPLMs and is of independent interest and importance. This theorem also
indicates that although there are several nonparametric functions and
their polynomial approximation deduces biases for the estimators of
each nonparametric component, these biases do not make the estimators
of $\beta$ biased under condition (C6) imposed on the number of knots.
\end{remark}

%s3 ###
\section{Focused information criterion and frequentist model averaging}{\label{SEC:model selection}}

In this section, based on the asymptotic result in Section \ref{SEC:framework and estimation}, we develop an FIC model selection
for GAPLMs, an FMA estimator, and propose a proper confidence interval
for the focus parameters.

%s3.1 ###
\subsection{Focused information criterion}\label{s31}

Let $\mu_{0}=\mu(\beta_0)=\mu(\betaczero,\delta/\sqrt{n})$ be a focus
parameter. Assume that the partial derivatives of $\mu(\beta_0)$ are
continuous in a neighborhood of $\betaczero$. Note that, in the
${\es}$th submodel, $\mu_{0}$ can be estimated by
$\widehat{\mu}_{\es}=\mu([I_{\dc},0_{\dc\times\du}]\Pi_{\es}^{\T}\widehat{\beta}_{\es},
[0_{\du\times\dc},I_{\du}]\Pi_{\es}^{\T}\widehat{\beta}_{\es})$. We now
show the asymptotic normality of  $\widehat{\mu}_{\es}$. Write
$\mathbf{R}_{\es}=\Pi_{\es}^{\T}(\Pi_{\es}\mathbf{D}\Pi_{\es}^{\T})^{-1}\Pi_{\es}$,
$\mu_{c}=\frac{\partial \mu (\beta_{ c},\betau)}{\partial \beta_{ c} }|_{\beta_{c}=\beta_{c, 0},\betau=0}$,
$\mu_{u}=\frac{\partial \mu (\beta_{c},\betau)}{\partial \beta_{u}}|_{\beta_{c}=\betaczero,\betau=0}$ and
$\mu_{\beta}=(\mu_{c}^{\T},\mu_{u}^{\T})^{\T}$.

\begin{theorem}\label{th2}
Under the local misspecification framework and conditions \textup{(C1)--(C11)} in the \hyperref[appe]{Appendix}, we have
\begin{eqnarray*}
\sqrt{n}(\widehat{\mu}_{\es}-\mu_{0})
&=&
-\mu_{\beta}^{\T}\mathbf{R}_{\es} \mathbf{G}_n+\mu_{\beta}^{\T} (\mathbf{R}_{\es} \mathbf{D}-I)
\pmatrix{
0\cr
\delta}+o_p(1)
\\
&\stackrel{d}\longrightarrow&
\Lambda_{\es}\equiv-\mu_{\beta}^{\T}\mathbf{R}_{\es}\mathbf{G}+\mu_{\beta}^{\T}(\mathbf{R}_{\es}\mathbf{D}-I)
\pmatrix{
0\cr
\delta}.
\end{eqnarray*}
\end{theorem}

Recall $\bfG\sim N(0, \Sigma)$. A direct calculation yields
%e3.1 ###
\begin{eqnarray}\label{eq:EGamma}
\qquad E(\Lambda_{\es}^2)
=
\mu_{\beta}^{\T}\left\{\matrix{\mathbf{R}_{\es}\bolds{\Sigma}\mathbf{R}_{\es}+(\mathbf{R}_{\es}\mathbf{D}-I)
\pmatrix{
0\cr
\delta}
\pmatrix{
0\cr
\delta}^{\T}(\mathbf{R}_{\es}\mathbf{D}-I)^{\T}}\right\}\mu_{\beta}.
\end{eqnarray}
Let $\widehat\delta$ be the estimator of $\delta$ by the full model.
Then, from Theorem \ref{th1}, we know that
\[
\widehat\delta=-[0,I]\mathbf{D}^{-1}\mathbf{G}_n+\delta+o_p(1).
\]
If we define $\Delta=-[0,I]\mathbf{D}^{-1}\mathbf{G}+\delta\sim N(\delta,[0,I]\mathbf{D}^{-1}\bolds{\Sigma}\mathbf{D}^{-1}[0,I]^{\T})$,
then $\widehat\delta\stackrel{d}\longrightarrow \Delta.$ Following
Claeskens and Hjort (\citeyear{ch2003}) and (\ref{eq:EGamma}), we define the FIC of
the ${\es}$th submodel as
\begin{eqnarray}\label{eq:FIC}
\hspace*{21pt}
\mathrm{FIC}_{\es}
&=&
\mu_{\beta}^{\T}\left\{
\mathbf{R}_{\es}\bolds{\Sigma}\mathbf{R}_{\es}+(\mathbf{R}_{\es}\mathbf{D}-I)
\pmatrix{
0 \cr
\widehat\delta}
\pmatrix{
0 \cr
\widehat\delta}^{\T}(\mathbf{R}_{\es}\mathbf{D}-I)^{\T}\right.\nonumber
\\[-8pt]\\[-8pt]
&&\hspace*{17pt}{}-
\left.(\mathbf{R}_{\es}\mathbf{D}-I)
\pmatrix{
0&0\cr
0&I_{\du}}\mathbf{D}^{-1}\bolds{\Sigma}\mathbf{D}^{-1}
\pmatrix{
0&0\cr
0&I_{\du}}
(\mathbf{R}_{\es}\mathbf{D}-I)^{\T}\right\}\mu_{\beta},\nonumber\hspace*{-21pt}
\end{eqnarray}
which is an approximately unbiased estimator of the mean squared error
when $\sqrt{n}\mu_{0}$ is estimated by $\sqrt{n}\widehat{\mu}_{\es}$.
This FIC can be used for choosing a proper submodel relying on the
parameter of interest.

%s3.2 ###
\subsection{Frequentist model averaging}\label{s32}

As mentioned previously, an average estimator is an alternative to a
model selection estimator. There are at least two advantages to the use
of an average estimator. First, an average estimator often reduces mean
square error in estimation because it avoids ignoring useful
information from the form of the relationship between response and
covariates and it provides a kind of insurance against selecting a very
poor submodel. Second, model averaging procedures can be more stable
than model selection, for which small changes in the data often lead to
a significant change in model choice. Similar discussions of this issue
appear in Bates and Granger (\citeyear{bg1969}) and Leung and Barron (\citeyear{lb2006}).

By choosing a submodel with the minimum value of FIC, the FIC
estimators of $\mu$ can be written as $
\widehat\mu_{\mathrm{FIC}}=\sum_{\es}\mathbf{I}$(FIC selects the $\es$th submodel)$\widehat{\mu}_{\es}$,
where $\mathbf{I}(\cdot)$, an indicator function, can be thought of
as a weight function depending on the data via $\widehat\delta$, %see (\ref{eq:FIC}), but
yet it just takes value either 0 or 1. To smooth estimators across
submodels, we may formulate the model average estimator of $\mu$
as
%e3.2 ###
\begin{equation}\label{eq9}
\widehat\mu=\sum_{\es}w(\es|\widehat\delta)\widehat{\mu}_{\es},
\end{equation}
where the weights $w(\es|\widehat\delta)$ take values in the interval
$[0,1]$ and their sum equals~1. It is readily seen  that smoothed AIC,
BIC and FIC estimators investigated in Hjort and Claeskens (\citeyear{hc2003}) and
Claeskens and Carroll (\citeyear{cc2007}) share this form. The following theorem
shows an asymptotic property for the general model average estimators
$\widehat\mu$ defined in (\ref{eq9}) under  certain conditions.

\begin{theorem}\label{th3}
Under the local misspecification framework and conditions
\textup{(C1)--(C11)} in the \hyperref[appe]{Appendix}, if the weight functions have at most
a countable number of discontinuities, then
\begin{eqnarray*}
\sqrt{n}(\widehat\mu-\mu_{0})
&=&
-\mu_{\beta}^{\T}\mathbf{D}^{-1}\mathbf{G}_n+\mu_{\beta}^{\T}\left\{Q(\widehat\delta)
\pmatrix{
0\cr
\widehat\delta}
-
\pmatrix{
0\cr
\widehat\delta}\right\}+o_p(1)
\\
&\stackrel{d}\longrightarrow&
\Lambda\equiv -\mu_{\beta}^{\T}\mathbf{D}^{-1}\mathbf{G}+\mu_{\beta}^{\T}\left\{Q(\Delta)
\pmatrix{
0\cr
\Delta}
-
\pmatrix{
0\cr
\Delta}\right\},
\end{eqnarray*}
where $Q(\cdot)=\sum_{\es}w(s|\cdot )\mathbf{R}_{\es}\mathbf{D}$ and
$\Delta$ is defined in Section \textup{\ref{s31}}.
\end{theorem}

Referring to the above theorems, we construct a confidence interval for
$\mu$ based on the model average estimator $\hat\mu$, as follows.
Assume that $\widehat\kappa^2$ is a consistent estimator of
$\mu_{\beta}^{\T}\mathbf{D}^{-1}\bolds{\Sigma}\mathbf{D}^{-1}\mu_{\beta}$.
It is easily seen that
\[
\left[\sqrt{n}(\widehat\mu-\mu_{0})-\mu_{\beta}^{\T}\left\{Q(\widehat\delta)
\pmatrix{
0\cr
\widehat\delta}
-
\pmatrix{
0\cr
\widehat\delta
}\right\}\right]\bigl/\widehat\kappa \stackrel{d}\longrightarrow N(0,1).
\]
If we define the lower bound ($\mathit{low}_n$) and upper bound ($\mathit{up}_n$) by
%e3.3 ###
\begin{equation}\label{eq:low of CI}
\widehat\mu-\mu_{\beta}^{\T}\left\{Q(\widehat\delta)
\pmatrix{
0\cr
\widehat\delta}
-
\pmatrix{
0\cr
\widehat\delta}\right\}\bigl/\sqrt{n}\mp z_{\jmath}\widehat\kappa/\sqrt{n},
\end{equation}
where $z_\jmath$ is the ${\jmath}$th standard
normal quantile, then we have $ \Pr\{ \mu_{0}\in (\mathit{low}_n,\break \mathit{up}_n)\}
\rightarrow 2\Phi(z_{\jmath})-1, $ where $\Phi(\cdot)$ is a standard
normal distribution function. Therefore, the interval $(\mathit{low}_n, \mathit{up}_n)$
can be used as a confidence interval for $\mu_0$ with asymptotic level
$2\Phi(z_{\jmath})-1$.

\begin{remark}\label{re:new}
Note that the limit distribution of $\sqrt{n}(\widehat\mu-\mu_{0})$ is
a nonlinear mixture of several normal variables. As argued in Hjort and
Claeskens (\citeyear{hc2006}), a direct construction of a confidence interval based
on Theorem \ref{th3} may not be easy. The confidence
interval based on (\ref{eq:low of CI}) %and (\ref{eq:up of CI})
is better in terms of coverage probability and computational
simplicity, as promoted in Hjort and Claeskens (\citeyear{hc2003}) and advocated by
Claeskens and Carroll (\citeyear{cc2007}).
\end{remark}

\begin{remark}\label{re:3}
A referee has asked whether the focus parameter can depend on the
nonparametric function $\eta_0$. Our answer is ``yes.'' For instance,
we consider a general focus parameter, $\eta_0(\bfx)+\mu_0$, a summand
of $\mu_0$, which we have studied, and a nonparametric value at $\bfx$.
We may continue to get an estimator of  $\eta_0(\bfx)+\mu_0$ by
minimizing (\ref{eq:FIC}) and then model-averaging estimators by
weighting the estimators of $\mu_0$ and $\eta_0$ as in (\ref{eq9}).
However, the underlying FMA estimators are not root-$n$ consistent
because the bias of these estimators is proportional to the bias of the
estimators of $\eta_0$, which is larger than $n^{-1/2}$, whereas we can
establish their rates of convergence using easier arguments than those
employed in the proof of Theorem~\ref{th3}. Even though the focus parameters
generally depend on $\mu_0$ and $\eta_0$ of form $H(\mu_0, \eta_0)$ for
a given function $H(\cdot, \cdot)$, the proposed method can be still applied.
However, to develop asymptotic properties for the corresponding FMA
estimators depends on the form of $H(\cdot, \cdot)$ and will require
further investigation. We omit the details. Our numerical studies below
follow these proposals when the focus parameters are related to the
nonparametric functions.
\end{remark}

%s4 ###
\section{Simulation study}\label{s4}

We generated $1000$ data sets consisting of $n=200$ and $400$
observations from the  GAPLM
\begin{eqnarray*}
\operatorname{logit}\{\Pr(Y_i=1)\}
&=&
\eta_{1}(\mathbf{X}_{i,1})+\eta_{2}(\mathbf{X}_{i,2})+\mathbf{Z}_i^{\T}\beta
\\
&=&
\sin(2\pi\mathbf{X}_{i,1})+5\mathbf{X}_{i,2}^4+3\mathbf{X}_{i,2}^2-2+\mathbf{Z}_i^{\T}\beta,\qquad i=1,\ldots,n,
\end{eqnarray*}
where: the true parameter
$\beta=\{1.5,2,{r_0}(2,1,3)/\sqrt{n}\}^{\T}$; $\mathbf{X}_{i,1}$ and
$\mathbf{X}_{i,2}$ are independently uniformly distributed on $[0,1]$;
$\mathbf{Z}_{i,1}, \ldots, \mathbf{Z}_{i,5}$ are normally distributed
with mean 0 and variance 1; when $\hbar_1\neq\hbar_2$, the correlation
between $\mathbf{Z}_{i,\hbar_1}$ and $\mathbf{Z}_{i,\hbar_2}$ is
$\varpi^{|\hbar_1-\hbar_2|}$ with $\varpi=0$ or $\varpi=0.5$;
$\mathbf{Z}_{i}$ is independent of $\mathbf{X}_{i,1}$ and
$\mathbf{X}_{i,2}$. We set the first two components of $\beta$ to be in
all submodels. The other three may or may not be present, so we have
$2^3=8$ submodels to be selected or averaged across. ${r_0}$ varies
from 1 or 4 to 7. Our focus parameters are (i)~$\mu_1=\beta_1$, (ii)
$\mu_2=\beta_2$, (iii)
$\mu_3=0.75\beta_1+0.05\beta_2-0.3\beta_3+0.1\beta_4-0.06\beta_5$ and
(iv)~$\mu_4=\eta_{1}(0.86)+\eta_{2}(0.53)+0.32\beta_1-0.87\beta_2-0.33\beta_3-0.15\beta_4+0.13\beta_5$.

The cubic B-splines have been used to approximate the two nonparametric
functions. We propose to select $J_{n}$ using a BIC procedure. Based on
condition (C6), the optimal order of $J_{n}$ can be found in the range
$(n^{1/(2\upsilon)}, n^{1/3})$. Thus, we propose to choose  the optimal
knot number, $J_{n}$, from a neighborhood of $n^{1/5.5}$. For our
numerical examples, we have used $[2/3N_{r},4/3N_{r}]$, where
$N_{r}=\operatorname{ceiling}(n^{1/5.5})$ and the function $\operatorname{ceiling}(\cdot)$ returns
the smallest integer not less than the corresponding element. Under the
full model, let the log-likelihood function be $l_{n}( N_{n})$. The
optimal knot number, $N^{\mathrm{opt}}_{n}$, is then the one which
minimizes the BIC value. That is,
%e4.1 ###
\begin{equation}\label{eqBIC}
N^{\mathrm{opt}}_{n}=\arg\min_{N_{n}\in[2/3N_{r},4/3N_{r}]}\{-2l_{n}(N_{n})+q_{n}\log n\},
\end{equation}
where $q_n$ is the total number of parameters.

Four model selection or model averaging methods are compared in this
simulation: AIC, BIC, FIC and the smoothed FIC (S-FIC). The smoothed
FIC weights we have used are
\[
w(\es|\hat\delta)=\exp\biggl(-\frac{\mathrm{FIC}_{\es}}{\mu_{\beta}^{\T}\bfD^{-1}\bf\Sigma\bfD^{-1}\mu_{\beta}}\biggr)
\Bigl/
\sum\limits_{\mathrm{all}\ \es'}\exp\biggl(-\frac{\mathrm{FIC}_{\es'}}{\mu_{\beta}^{\T}\bfD^{-1}\bf\Sigma\bfD^{-1}\mu_{\beta}}\biggr),
\]
a case of expression (5.4) in Hjort and Claeskens (\citeyear{hc2003}). When using
the FIC or S-FIC method, we estimate $\mathbf{D}^{-1}\bolds{\Sigma}
\mathbf{D}^{-1}$ by the covariance matrix of $\widehat\beta_{\mathrm{full}}$ and
estimate $\mathbf{D}$ by its sample mean, as advocated by Hjort and
Claeskens (\citeyear{hc2003}) and Claeskens and Carroll (\citeyear{cc2007}). Thus,
$\bolds{\Sigma}$ can be calculated straightforwardly. Note that the
subscript ``\textit{full}'' denotes the estimator using the full model.

In this simulation, one of our purposes is to see whether the
traditional selection methods like AIC and BIC lead to an overly
optimistic coverage probability (CP) of a claimed confidence interval
(CI). We consider a claimed 95\% confidence interval. The other purpose
is to check the accuracy of estimators in terms of their mean squared
errors (MSE)
${1}/{1000}\sum_{j}(\widehat\mu_a^{(j)}-\mu_a)^2$ for
$a=1,\ldots,4,$ where $j$ denotes the $j$th replication.
Our results are listed in Table \ref{ta:simulation}.

%t1 ###
\begin{table}[b]
\tabcolsep=0pt
\caption{Simulation results. Full: using all variables; CP: coverage
probability;  MSE: mean squared error}\label{ta:simulation}
{\fontsize{8pt}{10pt}\selectfont{
\begin{tabular*}{13.1cm}{@{\extracolsep{\fill}}lccd{1.2}cd{1.2}d{1.2}d{1.2}d{1.2}d{1.2}d{1.2}d{1.2}d{1.2}d{1.2}d{1.2}cd{1.2}d{1.2}d{1.2}@{}}
\hline\\[-15pt]
&&&\multicolumn{4}{c}{$\bolds{\mu_1}$}&\multicolumn{4}{c}{$\bolds{\mu_2}$}&\multicolumn{4}{c}{$\bolds{\mu_3}$}&\multicolumn{4}{c@{}}{$\bolds{\mu_4}$}\\[-5pt]
&&&\multicolumn{4}{c}{\hrulefill}&\multicolumn{4}{c}{\hrulefill}&\multicolumn{4}{c}{\hrulefill}&\multicolumn{4}{c@{}}{\hrulefill}\\
&&&\multicolumn{2}{c}{$\bolds{\varpi=0}$}&\multicolumn{2}{c}{$\bolds{\varpi=0.5}$}&\multicolumn{2}{c}{$\bolds{\varpi=0}$}&\multicolumn{2}{c}{$\bolds{\varpi=0.5}$}&\multicolumn{2}{c}{$\bolds{\varpi=0}$}&\multicolumn{2}{c}{$\bolds{\varpi=0.5}$}&\multicolumn{2}{c}{$\bolds{\varpi=0}$}&\multicolumn{2}{c@{}}{$\bolds{\varpi=0.5}$}\\[-5pt]
&&&\multicolumn{2}{c}{\hrulefill}&\multicolumn{2}{c}{\hrulefill}&\multicolumn{2}{c}{\hrulefill}&\multicolumn{2}{c}{\hrulefill}&\multicolumn{2}{c}{\hrulefill}&\multicolumn{2}{c}{\hrulefill}&\multicolumn{2}{c}{\hrulefill}&\multicolumn{2}{c@{}}{\hrulefill}\\
$\bolds{n}$&$\bolds{r_0}$&\multicolumn{1}{c}{\textbf{Method}}&\multicolumn{1}{c}{\textbf{CP}}&\multicolumn{1}{c}{\textbf{MSE}}&\multicolumn{1}{c}{\textbf{CP}}&\multicolumn{1}{c}{\textbf{MSE}}&\multicolumn{1}{c}{\textbf{CP}}&\multicolumn{1}{c}{\textbf{MSE}}&\multicolumn{1}{c}{\textbf{CP}}&\multicolumn{1}{c}{\textbf{MSE}}&\multicolumn{1}{c}{\textbf{CP}}&\multicolumn{1}{c}{\textbf{MSE}}&\multicolumn{1}{c}{\textbf{CP}}&\multicolumn{1}{c}{\textbf{MSE}}&\multicolumn{1}{c}{\textbf{CP}}&\multicolumn{1}{c}{\textbf{MSE}}&\multicolumn{1}{c}{\textbf{CP}}&\multicolumn{1}{c@{}}{\textbf{MSE}}\\
\hline
200&1&Full&0.9&0.33&0.9&0.49&0.89&0.49&0.88&0.8&0.9&0.22&0.9&0.32&0.92&2.25&0.91&2.92\\
&&AIC&0.91&0.31&0.9&0.46&0.9&0.45&0.87&0.76&0.89&0.21&0.88&0.3&0.91&2.15&0.9&2.77\\
&&BIC&0.92&0.28&0.9&0.4&0.91&0.39&0.88&0.71&0.9&0.19&0.88&0.26&0.92&1.98&0.9&2.66\\
&&FIC&0.92&0.28&0.93&0.39&0.92&0.33&0.91&0.79&0.92&0.19&0.92&0.25&&2&&2.66\\
&&S-FIC&0.93&0.28&0.93&0.41&0.93&0.4&0.92&0.68&0.93&0.19&0.92&0.26&&2&&2.61\\[3pt]
&4&Full&0.89&0.35&0.9&0.73&0.9&0.48&0.88&1.16&0.9&0.19&0.91&0.35&0.94&1.79&0.91&3.42\\
&&AIC&0.89&0.34&0.9&0.69&0.9&0.47&0.86&1.16&0.89&0.19&0.89&0.35&0.94&1.75&0.9&3.39\\
&&BIC&0.9&0.31&0.91&0.63&0.91&0.42&0.84&1.17&0.87&0.19&0.87&0.35&0.94&1.67&0.89&3.4\\
&&FIC&0.95&0.19&0.95&0.34&0.94&0.33&0.93&0.79&0.93&0.14&0.95&0.24&&1.52&&2.7\\
&&S-FIC&0.97&0.17&0.97&0.32&0.97&0.22&0.97&0.68&0.96&0.13&0.97&0.22&&1.32&&2.47\\[3pt]
&7&Full&0.89&0.46&0.9&1.02&0.89&0.66&0.87&2.04&0.9&0.2&0.92&0.41&0.92&2.26&0.92&5.32\\
&&AIC&0.89&0.46&0.9&1&0.89&0.65&0.86&2.04&0.9&0.2&0.91&0.41&0.92&2.24&0.91&5.28\\
&&BIC&0.89&0.44&0.91&0.93&0.9&0.62&0.86&1.92&0.9&0.2&0.88&0.41&0.92&2.18&0.91&4.79\\
&&FIC&0.94&0.21&0.97&0.36&0.94&0.33&0.95&0.79&0.95&0.12&0.97&0.19&&1.87&&2.98\\
&&S-FIC&0.97&0.12&0.98&0.22&0.97&0.16&0.98&0.63&0.98&0.09&0.98&0.15&&1.24&&2.57\\[3pt]
400&1&Full&0.93&0.07&0.92&0.1&0.93&0.11&0.93&0.15&0.93&0.05&0.92&0.07&0.94&0.52&0.94&0.67\\
&&AIC&0.94&0.07&0.92&0.1&0.93&0.1&0.91&0.14&0.93&0.04&0.91&0.07&0.94&0.51&0.93&0.66\\
&&BIC&0.94&0.06&0.93&0.09&0.94&0.09&0.91&0.14&0.93&0.04&0.91&0.06&0.94&0.5&0.93&0.65\\
&&FIC&0.94&0.06&0.93&0.09&0.94&0.09&0.94&0.15&0.94&0.04&0.93&0.06&&0.5&&0.65\\
&&S-FIC&0.95&0.06&0.93&0.09&0.94&0.1&0.93&0.14&0.94&0.04&0.94&0.06&&0.51&&0.64\\[3pt]
&4&Full&0.94&0.07&0.91&0.12&0.93&0.11&0.9&0.19&0.94&0.04&0.91&0.08&0.94&0.54&0.93&0.78\\
&&AIC&0.94&0.07&0.92&0.12&0.93&0.11&0.89&0.2&0.94&0.04&0.87&0.08&0.94&0.53&0.92&0.79\\
&&BIC&0.94&0.07&0.92&0.12&0.94&0.1&0.88&0.22&0.92&0.05&0.87&0.09&0.94&0.52&0.9&0.83\\
&&FIC&0.95&0.05&0.93&0.09&0.95&0.09&0.92&0.15&0.96&0.04&0.94&0.06&&0.49&&0.72\\
&&S-FIC&0.97&0.05&0.95&0.09&0.97&0.07&0.95&0.16&0.97&0.04&0.94&0.06&&0.46&&0.69\\[3pt]
&7&Full&0.92&0.08&0.9&0.14&0.93&0.11&0.91&0.21&0.94&0.04&0.91&0.08&0.94&0.52&0.93&0.82\\
&&AIC&0.92&0.08&0.9&0.14&0.92&0.11&0.91&0.21&0.94&0.04&0.89&0.08&0.94&0.51&0.93&0.81\\
&&BIC&0.93&0.08&0.91&0.13&0.93&0.11&0.89&0.22&0.94&0.04&0.86&0.09&0.94&0.5&0.92&0.82\\
&&FIC&0.94&0.06&0.92&0.1&0.93&0.09&0.93&0.15&0.94&0.04&0.93&0.07&&0.47&&0.68\\
&&S-FIC&0.95&0.05&0.96&0.07&0.95&0.06&0.96&0.12&0.96&0.03&0.96&0.05&&0.38&&0.6\\
\hline
\end{tabular*}}}
\end{table}

These results indicate that the performance of both the FIC and S-FIC,
especially the latter, is superior to that of AIC and BIC in terms of
CP and mean squared error (MSE), regardless of whether the focus
parameter depends on the nonparametric components or not. The CPs based
on FIC and S-FIC are generally close to the nominal level. When the
smallest CPs based on S-FIC and FIC are respectively 0.921 and 0.914,
the corresponding CPs of AIC and BIC are only 0.860 and 0.843,
respectively, which are much lower than the level 95\%. The CPs of both
S-FIC and FIC are higher than those from full models, but close to the
nominal level, whereas the intervals of FIC and S-FIC have the same
length as those from the full models because we estimate the unknown
quantities in
(\ref{eq:low of CI}) %and (\ref{eq:up of CI})
under the full model.

When ${r_0}$ gets bigger, the MSEs based on S-FIC are substantially
smaller than those obtained from other criteria. It is worth mentioning
that in Tables \ref{ta:simulation} and \ref{ta:over-fit}, we do not
report the CPs corresponding to FIC and S-FIC for $\mu_4$ because we do
not derive an asymptotic distribution for the proposed estimators of
this focus parameter.

As suggested by a referee, we now numerically examine the effects of
the number of knots on the performance of these criteria. We generalize
the data and conduct the simulation in the same way as above, but
oversmoothing and undersmoothing nonparametric terms by letting
$N_{r}=\operatorname{ceiling}(n^{1/3})$ and $N_{r}=\operatorname{ceiling}(n^{1/10})$, respectively.
The results corresponding to undersmoothing show a similar pattern as
in Table \ref{ta:simulation}. Note that  derivatives of all orders of
functions $\eta_1(X_{i,1})$ and $\eta_2(X_{i,2})$ exist and satisfy the
Lipschitz\vspace*{1pt} condition. $N_{r}=\operatorname{ceiling}(n^{1/10})$ is still in the range
$(n^{1/(2\upsilon)}, n^{1/3})$, so this similarity is not surprising
and supports our theory. However, oversmoothing of the nonparametric
functions causes significant changes and generally produces larger MSEs
but lower CPs, while all of the results show a preference for the S-FIC
and FIC. To save space, we report the results with $n=400$ and
${r_0}=4$ in Table \ref{ta:over-fit}, but omit other results, which
show similar features to those reported in Table \ref{ta:over-fit}.
\eject

%t2 ###
\begin{table}[b]
\tabcolsep=0pt
\tablewidth=\textheight
\tablewidth=\textwidth
\caption{Simulation results of overfitting with $n=400$ and ${r_0}=4$}\label{ta:over-fit}
\begin{tabular*}{\textwidth}{@{\extracolsep{\fill}}lcccccccc@{}}
\hline\\[-15pt]
&\multicolumn{4}{c}{$\bolds{\mu_1}$}&\multicolumn{4}{c@{}}{$\bolds{\mu_2}$}\\[-5pt]
&\multicolumn{4}{c}{\hrulefill}&\multicolumn{4}{c@{}}{\hrulefill}\\
&\multicolumn{2}{c}{$\bolds{\varpi=0}$}&\multicolumn{2}{c}{$\bolds{\varpi=0.5}$}&\multicolumn{2}{c}{$\bolds{\varpi=0}$}&\multicolumn{2}{c@{}}{$\bolds{\varpi=0.5}$}\\[-5pt]
&\multicolumn{2}{c}{\hrulefill}&\multicolumn{2}{c}{\hrulefill}&\multicolumn{2}{c}{\hrulefill}&\multicolumn{2}{c@{}}{\hrulefill}\\
\textbf{Method}&\textbf{CP}&\textbf{MSE}&\textbf{CP}&\textbf{MSE}&\textbf{CP}&\textbf{MSE}&\textbf{CP}&\textbf{MSE}\\
\hline
Full&0.864&0.131&0.852&0.232&0.852&0.211&0.840&0.365\\
AIC&0.869&0.129&0.863&0.226&0.851&0.207&0.805&0.381\\
BIC&0.884&0.117&0.872&0.210&0.863&0.186&0.770&0.409\\
FIC&0.942&0.086&0.917&0.154&0.922&0.131&0.874&0.300\\
S-FIC&0.952&0.081&0.932&0.149&0.946&0.123&0.916&0.300\\[6pt]
&\multicolumn{4}{c}{$\bolds{\mu_3}$}&\multicolumn{4}{c@{}}{$\bolds{\mu_4}$}\\[-5pt]
&\multicolumn{4}{c}{\hrulefill}&\multicolumn{4}{c@{}}{\hrulefill}\\
Full&0.884&0.073&0.863&0.138&0.928&1.055&0.910&1.548\\
AIC&0.874&0.073&0.813&0.142&0.931&1.053&0.909&1.571\\
BIC&0.863&0.077&0.782&0.152&0.929&1.028&0.897&1.606\\
FIC&0.914&0.060&0.915&0.107&&0.967&&1.443\\
S-FIC&0.949&0.064&0.921&0.110&&0.910&&1.361\\
\hline
\end{tabular*}
\end{table}

%s5 ###
\section{Real-world data analysis}\label{s5}

In this section, we apply our methods to a data set from a Pima Indian
diabetes study and perform some model selection and averaging
procedures. The data set is obtained from the UCI Repository of Machine
Learning Databases and selected from a larger data set held by the
National Institutes of Diabetes and Digestive and Kidney Diseases. The
patients under consideration are Pima Indian women at least 21 years
old and living near Phoenix, Arizona. The response variable, $Y$,
taking the value of 0 or 1, indicates a positive or negative test for
diabetes. The eight covariates are \textit{PGC} (plasma glucose
concentration after two hours in an oral glucose tolerance test),
\textit{DPF} (diabetes pedigree function), \textit{DBP} [diastolic
blood pressure (mm Hg)], \textit{NumPreg} (the number of times
pregnant), \textit{SI} [two-hour serum insulin (mu U$/$ml)],
\textit{TSFT} [triceps skin fold thickness (mm)], \textit{BMI}  (body
mass index [weight in kg$/$(height in m)$^2$]) and \textit{AGE} (years).
We then consider the following GAPLM for this data analysis:
\begin{eqnarray*}\label{model:PIW}
\operatorname{logit}\{\Pr(Y=1)\}
&=&
\eta_1(\mathit{BMI})+\eta_2(\mathit{AGE})+\beta_1\mathit{PGC}+\beta_2\mathit{DPF}
\\
&&{}+
\beta_3\mathit{DBP}+\beta_4\mathit{NumPreg}+\beta_5\mathit{SI}+\beta_6\mathit{TSFT},
\end{eqnarray*}
where \textit{AGE} and \textit{BMI} are set in nonparametric components and
the following Figure \ref{fig:etas} confirms that the effects of these
two covariates on the log odd are nonlinear. All covariates have been
centralized by sample mean and standardized by sample standard error.

%f1 ###
\begin{figure}[b]

\includegraphics{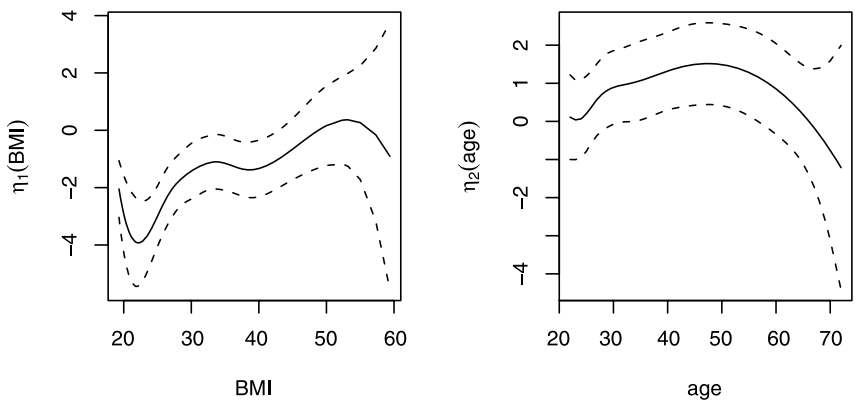}

\caption{The patterns of the nonparametric functions of
BMI and AGE (solid lines) with $\pm \mathrm{SE}$ (broken
lines).}\label{fig:etas}
\end{figure}

We first fit the model with all covariates using the polynomial spline
method introduced in Section \ref{SEC:framework and estimation}. The cubic B-splines have been used to
approximate the two nonparametric functions. The number of knots was
chosen using the BIC, presented in (\ref{eqBIC}). The fitted curves of
the two nonparametric components $\eta_1(\mathit{BMI})$ and $\eta_2(\mathit{AGE})$ are
depicted in Figure \ref{fig:etas}.  The estimated values of the
$\beta_i$'s, their standard error (SE) and corresponding z-values are
listed in Table \ref{ta:fitPIW}. The results indicate that \textit{PGC}
and \textit{DPF} are very significant, while the other four seem not to
be, so we run model selection and averaging on these four covariates.
Accordingly, there are $2^4=16$ submodels.

%t3 ###
\begin{table}
\tablewidth=227pt
\caption{Results for the diabetes study: estimated values, associated
standard errors and P-values obtained using the full model}\label{ta:fitPIW}
\begin{tabular}{@{\extracolsep{\fill}}ld{2.4}d{1.4}d{1.4}@{}}
\hline
&\multicolumn{1}{c}{\textbf{Estimated value}}&\multicolumn{1}{c}{\textbf{Standard error}}&\multicolumn{1}{c@{}}{$\bolds{P}$\textbf{-value}}\\
\hline
\textit{PGC}&1.1698&0.1236&0.0000\\
\textit{DPF}&0.3323&0.1029&0.0012\\
\textit{DBP}&-0.2662&0.1040&0.0104\\
\textit{NumPreg}&0.1887&0.1209&0.1184\\
\textit{SI}&-0.1511&0.1078&0.1610\\
\textit{TSFT}&0.0179&0.1135&0.8749\\
\hline
\end{tabular}
\end{table}

We now consider four focus parameters: $\mu_1=\beta_1$,
$\mu_2=\beta_2$,
$\mu_3=\eta_{1}(-1.501)+\eta_{2}(0.585)+0.028\beta_1-0.899\beta_2-1.570\beta_3+1.087\beta_4-0.223\beta_5-0.707\beta_6$
and
$\mu_4=\eta_{1}(-0.059)+\eta_{2}(1.363)+0.994\beta_1+0.423\beta_2+0.645\beta_3+1.117\beta_4-0.221\beta_5+0.055\beta_6$.
The first two are just the single coefficients of \textit{PGC} and
\textit{DPF}, the so-called two most significant linear components. The
second two are related to the nonparametric terms. Specifically
speaking, $\mu_3$ represents the log odd at $\mathit{BMI}=22.2$, the lowest
point of the estimated curve in the left panel of Figure~\ref{fig:etas}, and the corresponding means of other predictors when
$\mathit{BMI}=22.2$, while $\mu_4$ represents the log odd at $\mathit{AGE}=49$, the
highest point of the estimated curve in the right panel of Figure~\ref{fig:etas}, and the corresponding means of other predictors when
$\mathit{AGE}=49$. We label the potential 16 submodels ``0,'' ``3,''
``4,'' ``5,'' ``6,''$\dots,$ ``3456'' corresponding to a submodel which includes
(or not) \textit{DBP}, \textit{NumPreg}, \textit{SI} and \textit{TSFT}. The
results based on AIC, BIC and FIC methods are presented in Table
\ref{ta:focus parameters}. Regardless of focus parameter, the AIC and
BIC select submodels ``345'' and ``3,'' respectively. On the other hand,
the FIC prefers  submodels ``3,'' ``34,'' ``345'' and ``5'' when the focus
is on $\mu_1$, $\mu_2$, $\mu_3$ and $\mu_4$, respectively. It is
noticeable that submodel ``36'' is also competitive for $\mu_1$. We are
inclined to use submodel ``3'' since it has fewer parameters.

%
%t4 ###
\begin{sidewaystable}
\tabcolsep=0pt
\caption{Results for the diabetes study: AIC, BIC and FIC values, and estimators of focus parameters}\label{ta:focus parameters}
\begin{tabular*}{\textwidth}{@{\extracolsep{\fill}}ld{3.2}d{3.3}d{3.2}d{3.3}d{3.3}d{3.3}d{3.2}d{3.2}d{3.2}d{3.2}d{3.2}d{3.2}d{3.2}d{3.2}d{3.2}d{3.2}@{}}
\hline
&\multicolumn{1}{c}{\textbf{0}}&\multicolumn{1}{c}{\textbf{3}}&\multicolumn{1}{c}{\textbf{4}}&\multicolumn{1}{c}{\textbf{5}}&\multicolumn{1}{c}{\textbf{6}}&\multicolumn{1}{c}{\textbf{34}}&\multicolumn{1}{c}{\textbf{35}}&\multicolumn{1}{c}{\textbf{36}}&\multicolumn{1}{c}{\textbf{45}}&\multicolumn{1}{c}{\textbf{46}}&\multicolumn{1}{c}{\textbf{56}}&\multicolumn{1}{c}{\textbf{345}}&\multicolumn{1}{c}{\textbf{346}}&\multicolumn{1}{c}{\textbf{356}}&\multicolumn{1}{c}{\textbf{456}}&\multicolumn{1}{c@{}}{\textbf{3456}}\\
\hline
AIC&717.2&712.1&716.7&716.5&718.3&711.5&711.8&713.8&716.1&717.8&718.4&711.3^{\star}&713.2&713.7&718.0&713.3\\
BIC&791.3&790.7^{\star}&795.4&795.2&797.0&794.7&795.0&797.1&799.4&801.1&801.7&799.2&801.1&801.6&806.0&805.8\\[3pt]
 $\mu_1$-FIC&11.58&9.86^{\star}&13.69&11.07&11.4&10.97&11.74&9.86&11.63&13.41&11.36&11.16&10.96&12.17&12.08&11.54\\
 $\hat\mu_1$&1.09&1.11&1.09&1.15&1.092&1.11&1.17&1.11&1.15&1.09&1.15&1.17&1.11&1.17&1.15&1.17\\[3pt]
$\mu_2$-FIC&7.87&7.83&7.66&7.77&7.95&7.58^{\star}&7.77&8.07&7.79&7.93&7.97&7.80&8.00&7.98&8.01&7.99\\
$\hat\mu_2$&0.31&0.31&0.31&0.33&0.32&0.32&0.33&0.32&0.33&0.33&0.33&0.33&0.32&0.33&0.34&0.33\\[3pt]
$\mu_3$-FIC&261.5&51.9&143.7&245.8&219.1&38.0&48.7&47.9&144.7&122.5&235.1&37.7^{\star}&38.8&51.0&140.7&38.7\\
$\hat\mu_3$&-2.62&-2.23&-2.57&-2.70&-2.66&-2.17&-2.30&-2.26&-2.65&-2.61&-2.70&-2.24&-2.19&-2.29&-2.65&-2.23\\[3pt]
$\mu_4$-FIC&10.56&53.98&24.17&4.22^{\star}&9.82&30.70&30.08&51.82&35.28&23.98&6.02&31.38&30.71&30.10&35.43&31.93\\
$\hat\mu_4$&1.63&1.59&1.66&1.73&1.61&1.62&1.68&1.58&1.75&1.64&1.71&1.71&1.61&1.69&1.74&1.71\\
\hline
\end{tabular*}
$^{\star}$ denotes the minimal AIC, BIC or FIC values of the corresponding row.
\end{sidewaystable}

We further examine the predictive power of above model selection and
averaging methods through a cross-validation experiment. For each
patient in the data set, we use the AIC, BIC, FIC and S-FIC to carry
out estimations based on all of the other patients as a training
sample, and then predict the left-out observation. The prediction error
ratios (the ratio of the number of mistaken predictions to the sample
size) corresponding to AIC, BIC, FIC and S-FIC are 0.228, 0.225, 0.221
and 0.221, respectively. Both FIC and S-FIC show smaller prediction
errors than those of AIC and BIC, although the differences among these
errors are not substantial. These results indicate the superiority of
the FIC and S-FIC to the AIC and BIC.

%s6 ###
\section{Discussion}\label{s6}

We have proposed an effective procedure  using the polynomial spline
technique along with the model average principle to improve accuracy of
estimation in GAPLMs when uncertainty potentially appears. Our method
avoids any iterative algorithms and reduces computational challenges,
therefore its computational gain is remarkable. Most importantly, the
estimators of the linear components we have developed are still
asymptotically normal. Both theoretical and numerical studies show
promise for the proposed methods.

GAPLMs are generally enough to cover a variety of semiparametric models
such as partially linear additive models [Liang et al. (\citeyear{letal2008})] and
generalized partially linear models [Severini and  Staniswalis (\citeyear{ss1994})].
It is worth pointing out that GAPLMs do not involve any interaction
between nonparametric components (which may appear in a particular issue) and thus our current methods do not deal with this situation. We conjecture that
our procedure can be applied when the interactions may also be included
in the model search through tensor polynomial spline approximation, but
this extension poses additional challenges. How to develop model
selection and model averaging procedures in such a complex structure
warrants further investigation.

\begin{appendix}
\section*{Appendix}\label{appe}

Let $\|\cdot\|$ be the Euclidean norm and $\|\varphi\|_{\infty}=\sup_{m}|\varphi(m)|$ be the supremum norm of a function $\varphi$ on $[0,1] $. As in Carroll et al. (\citeyear{cetal1997}), we let
$q_{l}(m,y)=\frac{\partial^{l}Q\{g^{-1}(m),y\}}{\partial m^{l}}$, then
$q_{1}(m,y)={\partial Q\{ g^{-1}(m),y\}}/{\partial{m}}=\{y-g^{-1}(m)\}\rho_{1}(m)$ and $q_{2}(m,y)={\partial^{2}Q\{g^{-1}(m),y\}}/{\partial{m}^{2}}=\{y-g^{-1}(m)\}\rho_{1}^{\prime}(m)-\rho_{2}(m)$.

%s6.1 ###
\subsection{Conditions}\label{s61}
Let $r$ be a positive integer and $\nu\in(0,1]$ be such that
$\upsilon=r+\nu>1.5$. Let $\mathcal{H}$ be the collection of functions
$f$ on $[0,1]$ whose $r$th derivative, $f^{(r)}$,
exists and satisfies the Lipschitz condition of order $\nu$; that is,
\[
\bigl|f^{(r)}(m^{\ast})-f^{(r)}(m)\bigr|
\leq
C_1|m^{\ast}-m| ^{\nu }\qquad\mbox{for }0\leq m^{\ast},m\leq 1,
\]
where $C_1$ is a generic positive constant. In what follows,  $c$, $C$,
$c_{\bolds{\cdot}}$, $C_{\bolds{\cdot}}$ and $C_{\bolds{\cdot}}^{\ast}$ are all generic
positive constants. The following are the conditions needed to obtain
Theorems \ref{th1}--\ref{th3}:
\begin{enumerate}[(C10)]
\item[(C1)] each component function $\eta _{0,\alpha }\in\mathcal{H}$, $\alpha =1,\ldots,{p}$;
\item[(C2)] $q_{2}(m,y)<0$ and $c_{q}<|q_{2}(m,y)|<C_{q}$ for $m\in R$ and $y$ in the range of the response variable;
\item[(C3)] the function $\eta_{0}^{\prime\prime}(\cdot)$ is continuous;
\item[(C4)] the distribution of $\mathbf{X}$ is absolutely continuous and its density $f$ is bounded away from zero and infinity on $[0,1]^{{p}}$;
\item[(C5)] $E( \mathbf{ZZ}^{\T}|\mathbf{X}=\mathbf{x})$ exists and  $\mathbf{A}=E[\rho_{2}\{m_{0}(\mathbf{T})\}\mathbf{ZZ}^{\T}]$ is invertible, almost surely;
\item[(C6)] the number of interior knots $n^{1/(2\upsilon)}\ll J_{n}\ll n^{1/3}$;
\item[(C7)] $\lim_{n\rightarrow\infty}\frac{1}{n}\sum_{i=1}^n
\left({{\mathbf{B}(\mathbf{X}_{i})\mathbf{B}^{\T}(\mathbf{X}_{i})\enskip\mathbf{B}(\mathbf{X}_{i})\mathbf{Z}_{i}^{\T}}\atop
{\mathbf{Z}_{i}\mathbf{B}^{\T}(\mathbf{X}_{i})\enskip\mathbf{Z}_{i}\mathbf{Z}_{i}^{\T}}}\right)$ exists and is nonsingular;
\item[(C8)] for $\rho_{1}$ introduced in Section \ref{SEC:framework and estimation}, $|\rho_{1}(m_{0})|\leq C_{\rho}$ and
\[
|\rho_{1}(m)-\rho_{1}(m_{0})|\leq C_{\rho }^{*}|m-m_{0}|\qquad\mbox{for all }|m-m_{0}|\leq C_m;
\]
\item[(C9)] the matrix $\mathbf{D}$ is invertible almost surely;
\item[(C10)] the link function $g$ in model (\ref{model:GAPLM}) satisfies  $|\frac{d}{dm}g(m)|_{m=m_{0}}|\leq C_{g}$ and
\begin{eqnarray*}
&&\biggl|\frac{d}{dm}g^{-1}(m)-\frac{d}{dm}g^{-1}(m)\biggr| _{m=m_{0}}\biggl|
\\
&&\qquad\leq
C_{g}^{*}| m-m_{0}|\qquad\mbox{for all }|m-m_{0}|\leq C_m^{\ast};
\end{eqnarray*}
\item[(C11)] there exists a positive constant $C_{\varepsilon}$ such that $E(\varepsilon^{2}|\mathbf{T}=\bft)\leq C_{\varepsilon}$ almost
surely.
\end{enumerate}

%s6.2 ###
\subsection{Technical lemmas}\label{s62}

In the following, for any probability measure $P$, we define $L_{2}(P)=\{f\dvtx\int f^{2}\,dP<\infty\} $. Let $\mathcal{F}$ be a subclass
of $L_{2}(P)$. The bracketing number $\mathcal{N}_{[{}]}(\tau,\mathcal{F},L_{2}(P)) $ of $\mathcal{F}$ is defined as the smallest
value of $N$ for which there exist pairs of functions $\{[f_{j}^{L},f_{j}^{U}]\}_{j=1}^{N}$ with $\|f_{j}^{U}-f_{j}^{L}\|\leq\tau$, such that for each
$f\in\mathcal{F}$, there exists a $j\in\{ 1,\ldots,N\} $ such that $f_{j}^{L}\leq f\leq f_{j}^{U}$. Define the entropy integral
$J_{[{}]}(\tau,\mathcal{F},L_{2}(P))=\int_{0}^{\tau}\sqrt{1+\log\mathcal{N}_{[{}]}(\iota,\mathcal{F},L_{2}(P))}\,d\iota$.
Let $P_{n}$ be the empirical measure of $P$. Define $G_{n}= \sqrt{n}(
P_{n}-P) $ and $\| G_{n}\| _{\mathcal{F}}=\sup_{f\in \mathcal{F}}$ $|
G_{n}f| $ for any measurable class of functions $\mathcal{F}$.

We state or prove several preliminary lemmas first. Lemmas
\ref{LEM:Expection-entropy}--\ref{LEM:tildebeta} will be used to prove
the remaining lemmas. Lemmas
\ref{LEM:thetatilde-thetahat}--\ref{LEM:both} are used to prove Theorem~\ref{th1}. Theorems \ref{th2}--\ref{th3} are obtained from Theorem
\ref{th1}.

\begin{lemma}[{[Lemma 3.4.2 of van der Vaart and Wellner (\protect\citeyear{vdvw1996})]}]\label{LEM:Expection-entropy}
Let $M_{0}$ be a finite positive constant. Let $\mathcal{F}$ be a
uniformly bounded class of measurable functions such that
$Pf^{2}<\tau^{2}$ and $\|f\|_{\infty}<M_{0}$. Then
\[
E_{P}\|G_{n}\|_{\mathcal{F}}\leq C_{0}J_{[{}]}(\tau,\mathcal{F},L_{2}(P))
\biggl\{1+\frac{J_{[{}]}(\tau,\mathcal{F},L_{2}(P))}{\tau^{2}\sqrt{n}}M_{0}\biggr\},
\]
where $C_{0}$ is a finite constant not dependent on $n$.
\end{lemma}

\begin{lemma}[{[Lemma A.2 of Huang (\protect\citeyear{h1999})]}]\label{LEM:Entropy}
 For any $\tau>0$, let
$\Theta_{n}=\{\eta(\mathbf{x})+\mathbf{z}^{\T}\beta;\|\bolds{\beta}-\bolds{\beta}_{0}\|\leq\tau,\eta\in\mathcal{G}_{n},\|\eta-\eta_{0}\|_{2}\leq\tau\}$. Then, for any
$\iota\leq\tau$, $\log\mathcal{N}_{[]}(\iota,\Theta_{n},L_{2}(P))\leq c_0(J_{n}+\varrho)\log{\tau}/{\iota}$, where $c_{0}$ is a
finite constant not dependent on $n$.
\end{lemma}

Referring to the result of de Boor [(\citeyear{db2001}), page 149], for any function
$f\in\mathcal{H}$ and $n\geq 1$, there exists a function
$\widetilde{f}\in\mathcal{S}_{n}$ such that $\|\widetilde{f}-f\|_{\infty}\leq Ch^{\upsilon}$,
where $C$ is some fixed positive constant. From condition (C1), we can
find $\widetilde{\bolds{\gamma}}_{\es}=\{\widetilde{\gamma}_{\es,j,\alpha },j=-\varrho+1,\ldots,J_{n},\alpha=1,\ldots,p\}^{\T}$
and an additive spline function $\widetilde{\eta}_{\es}=\widetilde{\bolds{\gamma}}_{\es}^{\T}\mathbf{B}(\mathbf{x})\in\mathcal{G}_{n}$ such
that
%e6.1 ###
\begin{equation}\label{eq:tildeeta}
\|\widetilde{\eta}_{\es}-\eta_{0}\|_{\infty}=O(h^{\upsilon}) .
\end{equation}

Let $ \widetilde{\beta}_{\es}=\arg\max \frac{1}{n}\sum_{i=1}^{n}Q[
g^{-1}\{ \widetilde{\eta}_{\es}( \mathbf{X}_{i}) +(\Pi_{\es}\mathbf{Z}
_{i})^{\T}\beta_{\es}\} ,Y_{i}] $, $m_{0,i}= m_{0}( \mathbf{T}_{i})
=\eta _{0}( \mathbf{X}_{i}) +\mathbf{Z}_{i}^{\T}\bolds{\beta }_{0}$
and $\widetilde m_{\es,i}=\widetilde
m_{\es}(\mathbf{T}_i)=\widetilde\eta_{\es}(\mathbf{X}_i)+
\mathbf{Z}_i^{\T}\beta_0=\widetilde\gamma_{\es}^{\T}\bfB(\mathbf{X}_i)+\mathbf{Z}_i^{\T}\beta_0$.

\begin{lemma}\label{LEM:tildebeta}
Under the local misspecification framework and conditions \textup{(C1)--(C6)},
%e6.2 ###
\begin{equation}\label{eq:tildebeta}
\sqrt{n}\Pi_{\es}^{\T}(\widetilde{\beta}_{\es}-\Pi_{\es}\bolds{\beta}_{0})-\bar\Pi_{\es}^{\T}\bar\delta_{\es}
\stackrel{d}\longrightarrow
N(0,\mathbf{A}^{-1}\bolds{\Sigma}_{1}\mathbf{A}^{-1}),
\end{equation}
where $\bar\delta$ consists of the elements of $\delta$ that are not in
the ${\es}$th submodel, $\bar\pi_{\es}$ is the project matrix mapping
$\delta$ to $\bar\delta$, $\bar\Pi_{\es}=[0_{(d_{ u}-\duS)\times \dc},
\bar{\pi}_{\es}]$ and $\bolds{\Sigma }_{1}=E[q_{1}^{2}\{m_{0}(\mathbf{T})\}\mathbf{ZZ}^{\T}]$.
\end{lemma}

\begin{pf}
Let $\vartheta=\sqrt{n}\Pi_{\es}^{\T}(\beta_{\es}-\Pi_{\es}\bolds{\beta}_{0})-\bar\Pi_{\es}^{\T}\bar\delta_{\es}$ and
$\widetilde\vartheta=\sqrt{n}\Pi_{\es}^{\T}(\widetilde{\beta}_{\es}-\Pi_{\es}\bolds{\beta}_{0})-\bar\Pi_{\es}^{\T}\bar\delta_{\es}$. Note that
$\widetilde{\beta}_{\es}$ maximizes
%the quasi-likelihood function
$\frac{1}{n}\sum_{i=1}^{n}Q[g^{-1}\{\widetilde{\eta}_{\es}(\mathbf{X}_{i})+(\Pi_{\es}\mathbf{Z}_{i})^{\T}\beta_{\es}\},Y_{i}]$,
so $\widetilde\vartheta$ maximizes
\[
\widetilde{\ell}_{n}(\vartheta)
=
\sum_{i=1}^{n}[Q\{ g^{-1}(\widetilde{m}_{\es,i}+n^{-1/2}\vartheta^{\T}\mathbf{Z}_{i}),Y_{i}\}
-
Q\{g^{-1}(\widetilde{m}_{\es,i}),Y_{i}\}] .
\]
By Taylor expansion, one has\vspace*{-2pt}
$\widetilde{\ell}_{n}(\vartheta)=\frac{1}{\sqrt{n}}\sum_{i=1}^{n}q_{1}(\widetilde{m}_{\es,i},Y_{i})\vartheta^{\T}\mathbf{Z}_{i}+\frac{1}{2}\vartheta^{\T}\mathbf{A}_{n}\bolds{\vartheta}$, where
$\mathbf{A}_{n}=\frac{1}{n}\sum_{i=1}^{n}\{Y_{i}\rho_{1}^{\prime}(\widetilde{m}_{\es,i}+\zeta_{ni})-\rho_{3}(\widetilde{m}_{0i}+\zeta_{ni}^{\ast})\}\mathbf{Z}_{i}\mathbf{Z}_{i}^{\T}$
with $\zeta_{ni}$ and $\zeta_{ni}^{\ast}$ both lying between $0$ and $n^{-1/2}\vartheta^{\T}\mathbf{Z}_{i}$, and $\rho_{3}(m)=g^{-1}(m)\rho_{1}^{\prime}(m)-\rho_{2}(m)$.
From the proof of Theorem 2 in Carroll et al. (\citeyear{cetal1997}), $\mathbf{A}_{n}=-E[\rho_{2}\{m_{0}(\mathbf{T})\}\mathbf{ZZ}^{\T}]+o_{p}(1)=-\mathbf{A}+o_{p}(1)$ and
\begin{eqnarray*}
\frac{1}{\sqrt{n}}\sum_{i=1}^{n}q_{1}(\widetilde{m}_{\es,i},Y_{i})\mathbf{Z}_{i}
&=&
\frac{1}{\sqrt{n}}\sum_{i=1}^{n}q_{1}(m_{0,i},Y_{i})\mathbf{Z}_{i}
\\
&&{}+
\frac{1}{\sqrt{n}}\sum_{i=1}^{n}q_{2}(m_{0,i},Y_{i})\{\widetilde{\eta}_{\es}(\mathbf{X}_{i})-\eta_{0}(\mathbf{X}_{i})\}\mathbf{Z}_{i}
\\
&&{}+
O_{p}( n^{1/2}\|\widetilde{\eta}_{\es}-\eta_{0}\|_{\infty }^{2}).
\end{eqnarray*}
In addition, by (\ref{eq:tildeeta}) and conditions (C2), (C5) and (C6), we have
\begin{eqnarray*}
n^{-1/2}\sum_{i=1}^{n}q_{2}( m_{0,i},Y_i)\mathbf{Z}_{i}\{\widetilde{\eta}_\es(\mathbf{X}_{i})-\eta_{0}(\mathbf{X}_{i})\}=O_{p}(n^{1/2}h^{\upsilon})=o_{p}(1).
\end{eqnarray*}

Therefore, by the convexity lemma of Pollard (\citeyear{p1991}) and condition (C5),
one has $\widetilde{\vartheta}=\mathbf{A}^{-1}n^{-1/2}\sum_{i=1}^{n}q_{1}(m_{0,i},Y_i)\mathbf{Z}_{i}+o_{p}(1)$ and
$\operatorname{var}\{q_{1}\{ m_{0}(\mathbf{T}),Y \}\mathbf{Z}\}=E[q_{1}^{2}\{m_{0}(\mathbf{T}),Y \}\mathbf{ZZ}^{\T}]=\mathbf{\Sigma}_{1}$,
so (\ref{eq:tildebeta}) holds.
\end{pf}

Define $a_{n,h}=h^{\upsilon}+(n^{-1}\log n)^{1/2},$
$\theta_{\es}=(\gamma^{\T},\beta_{\es}^{\T})^{\T}$,
$\widetilde\theta_{\es}=(\widetilde\gamma_{\es}^{\T},\widetilde\beta_{\es}^{\T})^{\T}$ and
$\widehat\theta_{\es}=(\widehat\gamma_{\es}^{\T},\widehat\beta_{\es}^{\T})^{\T}$.

\begin{lemma}\label{LEM:thetatilde-thetahat}
Under the local misspecification framework and conditions \textup{(C1)--(C8)}, one has
$\|\widehat{\theta}_{\es}-\widetilde{\mathbf{\theta}}_{\es}\|=O_{p}(J_{n}^{1/2}a_{n,h}).
$
\end{lemma}

\begin{pf}
Note that
%e6.3 ###
\begin{equation}\label{eq:bartheta}
\frac{\partial\ell_{n}(\theta_{\es})}{\partial \theta_{\es}}\biggl|_{\theta_{\es}=\widehat\theta_{\es}}
-
\frac{\partial \ell_{n}(\theta_{\es})}{\partial \theta_{\es}}\biggl|_{\theta_{\es}=\widetilde\theta_{\es}}
=
\frac{\partial^{2}\ell_{n}(\theta_{\es})}{\partial\theta_{\es}\,\partial\theta_{\es}^{\T}}\biggl|_{\theta_{\es}=\bar\theta_{\es}}
(\widehat\theta_{\es}-\widetilde\theta_{\es}),
\end{equation}
with $\bar\theta_{\es}$ lying between $\widehat\theta_{\es}$ and
$\widetilde\theta_{\es}$. Recalling the equation
(\ref{DEF:quasilikelihood2}), one has
\[
\frac{\partial\ell_{n}(\theta_{\es})}{\partial\theta_{\es}}\biggl|_{\theta_{\es}=\widetilde\theta_{\es}}
=
\biggl\{\biggl(\frac{\partial\ell_{n}(\theta_{\es})}{\partial\gamma}\biggr)^{\T},
\biggl(\frac{\partial\ell_{n}(\theta_{\es})}{\partial\beta_{\es}}\biggr)^{\T}\biggr\}^{\T}\biggl|_{\theta_{\es}=\widetilde \theta_{\es}},
\]
where
\begin{eqnarray*}
\frac{\partial\ell_{n}(\theta_{\es})}{\partial\gamma}\biggl|_{\theta_{\es}=\widetilde\theta_{\es}}
&=&
\frac{1}{n}\sum_{i=1}^{n}q_{1}(m_{0,i},Y_{i})\mathbf{B}(\mathbf{X}_{i})
\\
&&{}+
\frac{1}{n}\sum_{i=1}^{n}q_{2}(\xi_{i},Y_{i})\{\widetilde\eta(\mathbf{X}_{i})-\eta_{0}(\mathbf{X}_{i})\}\mathbf{B}(\mathbf{X}_{i})
\\
&&{}+
\frac{1}{n}\sum_{i=1}^{n}q_{2}(\xi_{i},Y_{i})\bigl\{\Pi_{\es}^{\T}(\widetilde{\beta}_{\es}-\Pi_{\es}\beta_0)-\bar{\Pi}_{\es}^{\T}\bar{\delta}_{\es}/\sqrt{n}\bigr\}^{\T}\mathbf{Z}_{i}\mathbf{B}(\mathbf{X}_{i})
\end{eqnarray*}
and
\begin{eqnarray*}
\frac{\partial\ell_{n}(\theta_{\es})}{\partial\beta_{\es}}\biggl|_{\theta_{\es}=\widetilde\theta_{\es}}
&=&
\frac{1}{n}\sum_{i=1}^{n}q_{1}(m_{0,i},Y_{i})\Pi_{\es}\mathbf{Z}_{i}
\\
&&{}+
\frac{1}{n}\sum_{i=1}^{n}q_{2}(\xi^{\ast}_{i},Y_{i})\{\widetilde\eta(\mathbf{X}_{i})-\eta_{0}(\mathbf{X}_{i})\}\Pi_{\es}\mathbf{Z}_{i}
\\
&&{}+
\frac{1}{n}\sum_{i=1}^{n}q_{2}(\xi^{\ast}_{i},Y_{i})\bigl\{\Pi_{\es}^{\T}(\widetilde{\beta}_{\es}-\Pi_{\es}\beta_0)-\bar{\Pi}_{\es}^{\T}\bar{\delta}_{\es}/\sqrt{n}\bigr\}^{\T}\mathbf{Z}_{i}\Pi_{\es}\mathbf{Z}_{i},
\end{eqnarray*}
with $\xi_i$ and $\xi^{\ast}_i$ both lying between $m_{0,i}$ and
$\widetilde m_{\es,i}$.
According to the Bernstein inequality and condition (C8),
\begin{eqnarray*}
\Biggl\|\frac{1}{n}\sum_{i=1}^{n}q_{1}(m_{0,i},Y_{i})\mathbf{B}(\mathbf{X}_{i})\Biggr\|_{\infty }
&=&
\max_{-\varrho+1\leq j\leq J,1\leq\alpha\leq p}\frac{1}{n}\Biggl|\sum_{i=1}^{n}\rho_{1}( m_{0,i})B_{j,\alpha}(X_{i\alpha})\varepsilon_{i}\Biggr|
\\
&=&
O_{p}\{( n^{-1}\log n)^{1/2}\}.
\end{eqnarray*}
And, by (\ref{eq:tildeeta}), Lemma \ref{LEM:tildebeta} and condition
(C2), one has
\[
\frac{1}{n}\sum_{i=1}^{n}\|q_{2}(\mathbf{\xi}_{i},Y_{i})\{\widetilde\eta_{\es}(\mathbf{X}_{i})-\eta_{0}(\mathbf{X}_{i})\}\mathbf{B}(\mathbf{X}_{i})\|_{\infty}=O_{p}(h^{\upsilon})
\]
and
\[
\frac{1}{n}\sum_{i=1}^{n}\bigl\|q_{2}(\mathbf{\xi}_{i},Y_{i})
\bigl\{\Pi_{\es}^{\T}(\widetilde{\beta}_{\es}-\Pi_{\es}\beta_0)-\bar{\Pi}_{\es}^{\T}\bar{\delta}_{\es}/\sqrt{n}\bigr\}^{\T}\mathbf{Z}_{i}\mathbf{B}(\mathbf{X}_{i})\bigr\|_{\infty}=O_{p}(n^{-1/2}).
\]
Therefore, $\|{\partial\ell_{n}(\theta_{\es})}/{\partial\gamma}|_{\theta_{\es}=\widetilde\theta_{\es}}\|_{\infty}=O_{p}(a_{n,h})$. Similarly, we can prove
\[
\biggl\|\frac{\partial\ell_{n}(\theta_{\es})}{\partial\beta_{\es}}\biggl|_{\theta_{\es}=\widetilde\theta_{\es}}\biggr\|_{\infty}
=
O_{p}\bigl(h^{\upsilon}+(n^{-1}\log n)^{1/2}\bigr).
\]
Thus,
%e6.4 ###
\begin{equation}\label{eq:diftheta}
\biggl\|\frac{\partial\ell_{n}(\theta_{\es})}{\partial\theta_{\es}}\biggl|_{\theta_{\es}=\widetilde\theta_{\es}}\biggr\|_{\infty}
=
O_{p}(a_{n,h}).
\end{equation}

Let $\bar m_{\es,i}=\bar m_{\es}(\mathbf{T}_i)=\bar\theta^{\T}(\mathbf{B}^{\T}(\mathbf{X}_i),(\Pi_{\es}\mathbf{Z}_i)^{\T})^{\T}$.
For the second order derivative, one has
\begin{eqnarray*}
\frac{\partial ^{2}\ell_{n}(\theta_{\es})}{\partial \theta_{\es}\,\partial\theta_{\es}^{\T}}\biggl|_{\theta_{\es}=\bar\theta_{\es}}
&=&
\left.
\pmatrix{
\displaystyle\frac{\partial^{2}\ell_{n}(\mathbf{\theta}_{\es})}{\partial\mathbf{\gamma}\,\partial\mathbf{\gamma}^{\T}}&\displaystyle\frac{\partial^{2}\ell_{n}(\mathbf{\theta}_{\es})}{\partial\mathbf{\gamma}\,\partial\mathbf{\beta}_{\es}^{\T}}\vspace*{2pt}\cr
\displaystyle\frac{\partial^{2}\ell_{n}(\mathbf{\theta}_{\es})}{\partial\mathbf{\beta}_{\es}\,\partial\mathbf{\gamma}^{\T}}&\displaystyle\frac{\partial^{2}\ell_{n}(\mathbf{\theta}_{\es})}{\partial\mathbf{\beta}_{\es}\,\partial\mathbf{\beta}_{\es}^{\T}}
}\right|_{\theta_{\es}=\bar\theta_{\es}}
\\
&=&
\frac{1}{n}\sum_{i=1}^{n}q_{2}( \bar{m}_{\es,i},Y_{i})
\left\{\pmatrix{
\mathbf{B}(\mathbf{X}_{i})\mathbf{B}^{\T}(\mathbf{X}_{i})&\mathbf{B}(\mathbf{X}_{i})\mathbf{Z}_{i}^{\T}\vspace*{2pt}\cr
\mathbf{Z}_{i}\mathbf{B}^{\T}(\mathbf{X}_{i})&\mathbf{Z}_{i}\mathbf{Z}_{i}^{\T}
}\right\},
\end{eqnarray*}
by which, along with conditions (C2) and (C7), we know that the matrix
$\frac{\partial^{2}\ell_{n}(\theta_{\es})}{\partial\theta_{\es}\,\partial\theta_{\es}^{\T}}|_{\theta_{\es}=\bar\theta_{\es}}$ is nonsingular in probability.
So, according to (\ref{eq:bartheta}) and (\ref{eq:diftheta}), we have
completed the proof.
\end{pf}

Define $ \mathcal{M}_{n}=\{m(\mathbf{x},\mathbf{z})=\eta(\mathbf{x})+\mathbf{z}^{\T}\beta\dvtx\eta\in\mathcal{G}_{n}\}$ and
a class of functions $\mathcal{A}(\tau)=\{\rho_{1}(m(\mathbf{t}))\mathbf{\psi}(\mathbf{t})\dvtx m\in\mathcal{M}_{n},\|m-m_{0}\|\leq\tau\}$.

\begin{lemma}\label{LEM:both}
Under the local misspecification framework and conditions
\textup{(C1)--(C8)}, we have
%e6.6 ###
%e6.5 ###
\begin{eqnarray}
\label{eq:both1}\frac{1}{n}\sum_{i=1}^{n}\{\widehat{\eta}_{\es}(\mathbf{X}_{i})-\eta _{0}(\mathbf{X}_{i})\}\rho _{1}( m_{0,i})\mathbf{\psi}(\mathbf{T}_{i})
&=&
o_{p}( n^{-1/2}),
\\
\label{eq:both2}\frac{1}{n}\sum_{i=1}^{n}\rho_{1}( m_{0,i})\mathbf{\psi}(\mathbf{T}_{i})\mathbf{\Gamma}(\mathbf{X}_{i})^{\T}\Pi_{\es}^{\T}(\widehat{\beta}_{\es}-\Pi_{\es}\beta_0)
&=&
o_{p}(n^{-1/2}) .
\end{eqnarray}
\end{lemma}

\begin{pf}
Noting that $\mathbf{\psi}$ and $\rho_1$ are fixed bounded functions
under condition~(C8), by Lemma \ref{LEM:Entropy}, similar to the proof
of Corollary A.1 in Huang (\citeyear{h1999}), we can show, for any
$\iota\leq\tau$,
$\log\mathcal{N}_{[]}(\iota,\mathcal{A}(\tau),\|\cdot\|)\leq
c_0((J_{n}+\varrho)\log({\tau}/{\iota})+\log(\iota^{-1}))$, so the
corresponding entropy integral satisfies
$J_{[]}(\tau,\mathcal{A}(\tau),\|\cdot\|)\leq c_0\tau
\{(J_{n}+\varrho)^{1/2}+(\log\tau^{-1})^{1/2}\}$. According to Lemma
\ref{LEM:thetatilde-thetahat},
$\|\widehat\eta_{\es}-\widetilde\eta_{\es}\|_{2}^{2}=(\widehat\gamma_{\es}-\widetilde\gamma_{\es})^{\T}\sum_{i=1}^{n}E\{\mathbf{B}(\mathbf{X}_{i})\mathbf{B}^{\T}(\mathbf{X}_{i})\}
(\widehat\gamma_{\es}-\widetilde\gamma_{\es})/n\leq
C_7\|\widehat\gamma_{\es}-\widetilde\gamma_{\es}\|_{2}^{2}$, thus
$\|\widehat\eta_{\es}-\widetilde\eta_{\es}\|_{2}=O_{p}(J_{n}^{1/2}a_{n,h})$
and
$\|\widehat\eta_{\es}-\eta_{0}\|_{2}\leq\|\widehat\eta_{\es}-\widetilde\eta_{\es}\|_{2}+\|\widetilde\eta_{\es}-\eta_{0}\|_{2}=O_{p}(J_{n}^{1/2}a_{n,h})$.
Now, by Lemma 7 of Stone (\citeyear{s1986}),
%e6.7 ###
\begin{equation}\label{eq:etainfty}
\|\widehat\eta_{\es}-\eta_{0}\|_{\infty}\leq C_8J_n^{1/2}\|\widehat\eta_{\es}-\eta_{0}\|_{2}=O_{p}(J_{n}a_{n,h}).
\end{equation}
Thus, by Lemma \ref{LEM:Expection-entropy}, together with conditions
(C1) and (C6), we have
\begin{eqnarray*}
&&E\Biggl|\frac{1}{n}\sum_{i=1}^{n}\{\widehat\eta_{\es}(\mathbf{X}_{i})-\eta_{0}(\mathbf{X}_{i})\}\rho_{1}(m_{0,i})\mathbf{\psi}(\mathbf{T}_{i})
\\
&&\hspace*{17pt}{}-
E[\{\widehat\eta_{\es}(\mathbf{X})-\eta_{0}(\mathbf{X})\}\rho_{1}\{m_{0,i}\mathbf{\psi}(\mathbf{T})\}]\Biggr|=o(n^{-1/2}).
\end{eqnarray*}
In addition, by the definition of $\mathbf{\psi}$,
$E[\phi(\mathbf{X})\rho_{1}\{m_{0}(\mathbf{T})\}\mathbf{\psi}(\mathbf{T})]=0$ for any measurable function $\phi$. Hence
(\ref{eq:both1}) holds. Similarly, (\ref{eq:both2}) follows from Lemmas
\ref{LEM:Expection-entropy}--\ref{LEM:thetatilde-thetahat}.
\end{pf}

%s6.3 ###
\subsection{\texorpdfstring{Proof of Theorem \protect\ref{th1}}{Proof of Theorem 1}}
\label{s63}
Let $\widehat m_{\es,i}=\widehat m_{\es}(\mathbf{T}_i)=\widehat\eta_{\es}(\mathbf{X}_i)+\widehat\beta_{\es}^{\T}\Pi_{\es}\mathbf{Z}_i$. For any
$\mathbf{v}\in R^{\dc+\duS}$, define
$\widehat{m}_{\es}(\mathbf{v})=\widehat{m}_{\es}(\mathbf{x},\Pi_{\es}\mathbf{z})+
\mathbf{v}^{\T}\{\Pi_{\es}\mathbf{z}-\Pi_{\es}\mathbf{\Gamma}(\mathbf{x})\}=\break
\widehat{m}_{\es}(\mathbf{x},\Pi_{\es}\mathbf{z})+\mathbf{v}^{\T}\Pi_{\es}\psi(\mathbf{t})$.
Note that when $\mathbf{v}=0$, $\widehat{m}_{\es}(\mathbf{v})$
maximizes\break ${1}/{n}\sum_{i=1}^nQ[g^{-1}\{m_{\es}(\mathbf{T}_i)\}, Y_i]$ for all
$m_{\es}\in\{m_{\es}(\mathbf{x},\mathbf{z})=\eta(\mathbf{x})+(\Pi_{\es}\mathbf{z})^{\T}\beta_{\es}\dvtx\eta\in\mathcal{G}_{n}\}$, by which
\begin{eqnarray}\label{eq:123}
\hspace*{25pt}0
&=&
\frac{\partial}{\partial\mathbf{v}}\ell_{n}(\widehat{m}_{\es}(\mathbf{v}))\biggl|_{\mathbf{v}=0}\nonumber
\\
&=&
\frac{1}{n}\sum_{i=1}^{n}\{Y_{i}-g^{-1}(\widehat{m}_{\es,i})\}\rho_{1}(\widehat{m}_{\es,i})\Pi_{\es}\mathbf{\psi}(\mathbf{T}_{i}).\nonumber
\\
&=&
\frac{1}{n}\sum_{i=1}^{n}q_1(m_{0,i},Y_i)\Pi_{\es}\mathbf{\psi}(\mathbf{T}_{i})+\frac{1}{n}\sum_{i=1}^{n}\varepsilon_{i}\{\rho_{1}(\widehat{m}_{\es,i})-\rho_{1}(m_{0,i})\}\Pi_{\es}\mathbf{\psi}(\mathbf{T}_{i})
\\
&&{}-
\frac{1}{n}\sum_{i=1}^{n}\{g^{-1}(\widehat{m}_{\es,i})-g^{-1}( m_{0,i})\}\rho_{1}(\widehat{m}_{\es,i})\Pi_{\es}\mathbf{\psi}(\mathbf{T}_{i})\nonumber
\\
&\equiv&
I+\mathit{II}-\mathit{III}\nonumber.
\end{eqnarray}
Note that for the second term $E[\varepsilon_{i}\{\rho_{1}(\widehat{m}_{\es,i})-\rho_{1}(m_{0,i})\}\Pi_{\es}\mathbf{\psi}(\mathbf{T}_{i})] =0$.
From Lemma \ref{LEM:tildebeta}, (\ref{eq:both2}) and (\ref{eq:etainfty}), we have $\|\widehat{m}_{\es}-m_{0}\|_{\infty}=O_{p}(J_{n}^{1/2}a_{n,h})$, so,
by condition (C8), $\|\rho_{1}(\widehat{m}_{\es})-\rho_{1}(m_{0})\|_{\infty}=O_{p}(J_{n}^{1/2}a_{n,h}). $ Now, by the Bernstein inequality, under
condition (C11), we show that
%e6.8 ###
\begin{equation}\label{eq:II}
\mathit{II}=\frac{1}{n}\sum_{i=1}^{n}\varepsilon_{i}\{\rho_{1}(\widehat{m}_{\es,i})-\rho_{1}(m_{0,i})\}\Pi_{\es}\mathbf{\psi}(\mathbf{T}_{i})
=
o_{p}( n^{-1/2}).
\end{equation}
Express the third term as\vspace*{-1pt}
\begin{eqnarray*}
\mathit{III}
&=&
\frac{1}{n}\sum_{i=1}^{n}\{g^{-1}(\widehat{m}_{\es,i})-g^{-1}(m_{0,i})\}\rho_{1}(\widehat{m}_{\es,i})\Pi_{\es}\mathbf{\psi}(\mathbf{T}_{i})\vspace*{-2pt}
\\
&=&
\frac{1}{n}\sum_{i=1}^{n}(\widehat{m}_{\es,i}-m_{0,i})\rho_{1}(m_{0,i})\Pi_{\es}\mathbf{\psi}(\mathbf{T}_{i})\vspace*{-2pt}
\\
&&{}+
\frac{1}{n}\sum_{i=1}^{n}\{g^{-1}(\widehat{m}_{\es,i})-g^{-1}(m_{0,i})-(\widehat{m}_{\es,i}-m_{0,i})\}\rho_{1}(m_{0,i})\Pi_{\es}\mathbf{\psi}(\mathbf{T}_{i})\vspace*{-2pt}
\\
&&{}+
\frac{1}{n}\sum_{i=1}^{n}\{g^{-1}(\widehat{m}_{\es,i})-g^{-1}(m_{0,i})\}\{\rho_{1}(\widehat{m}_{\es,i})-\rho_{1}( m_{0,i})\}\Pi_{\es}\mathbf{\psi}(\mathbf{T}_{i})\vspace*{-2pt}
\\
&\equiv&
\mathit{III}_{1}+\mathit{III}_{2}+\mathit{III}_{3}.\vspace*{-1pt}
\end{eqnarray*}
From Lemma \ref{LEM:both}, a direct simplification yields\vspace*{-1pt}
\begin{eqnarray*}
\mathit{III}_{1}
&=&
\frac{1}{n}\sum_{i=1}^{n}\{\widehat{\eta}_{\es}(\mathbf{X}_i)+\widehat{\beta}_{\es}^{\T}\Pi_{\es}\mathbf{Z}_i-\eta_0(\mathbf{X}_i)-\beta_0^{\T}\mathbf{Z}_i\}\rho_1(m_{0,i})\Pi_{\es}\psi(\mathbf{T}_i)\vspace*{-2pt}
\\
&=&
\frac{1}{n}\sum_{i=1}^{n}\{\widehat{\eta}_{\es}(\mathbf{X}_i)-\eta_0(\mathbf{X}_i)+(\Pi_{\es}^{\T}\widehat{\beta}_{\es}-\beta_0)^{\T}\psi(\mathbf{T}_i)\vspace*{-2pt}
\\
&&\hspace*{103pt}{}+
(\Pi_{\es}^{\T}\widehat{\beta}_{\es}-\beta_0)^{\T}\mathrm{\Gamma}(\mathbf{X}_i)\}\rho_1(m_{0,i})\Pi_{\es}\psi(\mathbf{T}_i)\vspace*{-2pt}
\\
&=&
\frac{1}{n}\sum_{i=1}^{n}\{\widehat{\eta}_{\es}(\mathbf{X}_i)-\eta_0(\mathbf{X}_i)\}\rho_1(m_{0,i})\Pi_{\es}\psi(\mathbf{T}_i)\vspace*{-2pt}
\\
&&{}+
\frac{1}{n}\sum_{i=1}^{n}\rho_1(m_{0,i})\Pi_{\es}\psi(\mathbf{T}_i)\psi(\mathbf{T}_i)^{\T}\left\{\Pi_{\es}^{\T}\widehat\beta_{\es}-
\pmatrix{
\betaczero\cr
0}\right\}\vspace*{-2pt}
\\
&&{}-
\frac{1}{n}\sum_{i=1}^{n}\rho_1(m_{0,i})\Pi_{\es}\psi(\mathbf{T}_i)\psi(\mathbf{T}_i)^{\T}[0,I]^{\T}\delta/\sqrt{n}\vspace*{-2pt}
\\
&&{}+
\frac{1}{n}\sum_{i=1}^{n}\rho_1(m_{0,i})\Pi_{\es}\psi(\mathbf{T}_i)\mathrm{\Gamma}(\mathbf{X}_i)^{\T}\Pi_{\es}^{\T}(\widehat{\beta}_{\es}-\Pi_{\es}\beta_0)\vspace*{-2pt}
\\
&&{}+
\frac{1}{n}\sum_{i=1}^{n}\rho_1(m_{0,i})\Pi_{\es}\psi(\mathbf{T}_i)\mathrm{\Gamma}(\mathbf{X}_i)^{\T}\bar\Pi_{\es}^{\T}\bigl(-\bar{\delta}_{\es}/\sqrt{n}\bigr)\vspace*{-2pt}
\\
&=&
\frac{1}{n}\sum_{i=1}^{n}\rho_1(m_{0,i})\Pi_{\es}\psi(\mathbf{T}_i)\psi(\mathbf{T}_i)^{\T}\Pi_{\es}^{\T}\left\{\widehat{\beta}_{\es}-
\pmatrix{
\betaczero\cr
0
}\right\}\vspace*{-2pt}
\\
&&{}-
\frac{1}{\sqrt{n}}\frac{1}{n}\sum_{i=1}^{n}\rho_1(m_{0,i})\Pi_{\es}\psi(\mathbf{T}_i)\psi(\mathbf{T}_i)^{\T}[0,I]^{\T}\delta+o_p(n^{-1/2}).
\end{eqnarray*}
In addition, from conditions (C8) and (C10), referring to the proof of
(\ref{eq:both1}), we have $\mathit{III}_{2}=o_p(n^{-1/2})$ and
$\mathit{III}_{3}=o_p(n^{-1/2})$. Therefore,
\begin{eqnarray}\label{eq:III}
\qquad\qquad\mathit{III}
&=&
[E\{\rho_1(m_{0})\Pi_{\es}\psi(\mathbf{T})\psi(\mathbf{T})^{\T}\Pi_{\es}^{\T}\}+o_p(1)]\left\{\widehat{\beta}_{\es}-
\pmatrix{
\betaczero\cr
0
}\right\}\nonumber
\\[-8pt]\\[-8pt]
&&{}-
\frac{1}{\sqrt{n}}[E\{\rho_1(m_{0})\Pi_{\es}\psi(\mathbf{T})\psi(\mathbf{T})^{\T}[0,I]^{\T}\}+o_p(1)]\delta+o_p(n^{-1/2}).\nonumber
\end{eqnarray}
Thus, by\vspace*{0.5pt} combining (\ref{eq:123}), (\ref{eq:II}), (\ref{eq:III}) and
condition (C9), the desired distribution of $\widehat\beta_{\es}$
follows.\vspace*{-2pt}

%s6.4 ###
\subsection{\texorpdfstring{Proof of Theorem \protect\ref{th2}}{Proof of Theorem 2}}
\label{s64}

By the Taylor expansion,
$\mu_{0}=\mu(\betaczero,\delta/\sqrt{n})=\mu(\betaczero,0)+\mu_u^{\T}\delta/\sqrt{n}+o(n^{-1/2})$
and
\begin{eqnarray*}
\widehat{\mu}_{\es}
&=&
\mu([I,0]\Pi_{\es}^{\T}\widehat{\beta}_{\es},
[0,I]\Pi_{\es}^{\T}\widehat{\beta}_{\es})
\\
&=&
\mu(\betaczero,0)+\mu_{\beta}^{\T}\left\{\Pi_{\es}^{\T}\widehat{\beta}_{\es}-
\pmatrix{
\betaczero\cr
0
}\right\}+o_p(n^{-1/2}),
\end{eqnarray*}
where the second equation follows from the asymptotic normality of
$\widehat{\beta}_{\es}$. Thus, by Theorem \ref{th1},
\begin{eqnarray*}
\sqrt{n}(\widehat{\mu}_{\es}-\mu_{0})&=&\mu_{\beta}^{\T}\left\{\Pi_{\es}^{\T}\widehat{\beta}_{\es}-
\pmatrix{
\betaczero\cr
0
}\right\}-\mu_u^{\T}\delta + o_p(1)
\\
%0_{\duS\times 1}
%)\}-\mu_u^{\T}\delta + o_p(1)\\
&=&
-\mu_{\beta}^{\T}\mathbf{R}_{\es} \mathbf{G}_n+\mu_{\beta}^{\T}\mathbf{R}_{\es}\mathbf{D}
\pmatrix{
0\cr
\delta}-\mu_u^{\T}\delta+o_p(1)
\\
&\stackrel{d}\longrightarrow&
-\mu_{\beta}^{\T}\mathbf{R}_{\es}\mathbf{G} +\mu_{\beta}^{\T}(\mathbf{R}_{\es}\mathbf{D}-I)
\pmatrix{
0\cr
\delta}.
\end{eqnarray*}
Thus, the proof is complete.\vspace*{-3pt}

%s6.5 ###
\subsection{\texorpdfstring{Proof of Theorem \protect\ref{th3}}{Proof of Theorem 3}}
\label{s65}

Recalling the definitions of $\Pi_{\es}$ and $\mathbf{R}_{\es}$, we
have
\[
\mathbf{R}_{\es}\mathbf{D}
\pmatrix{
I&0\cr
0&0_{\du\times \du}}
=
\mathbf{R}_{\es}\mathbf{D}\Pi_{\es}^{\T}
\pmatrix{I&0\cr
0&0_{\duS\times \du}}
=
\pmatrix{
I&0\cr
0&0_{\du\times \du}},
\]
which, along with the definition of $\widehat\delta$  and Theorem \ref{th2}, indicates that
\begin{eqnarray*}
\sqrt{n}(\widehat\mu-\mu_{0})
&=&
\sum_{\es}w(\es|\widehat\delta)\sqrt{n}(\widehat{\mu}_{\es}-\mu_{0})
\\
&=&
\sum_{\es}w(\es|\widehat\delta)\left\{-\mu_{\beta}^{\T}\mathbf{R}_{\es}\mathbf{G}_n+\mu_{\beta}^{\T}\mathbf{R}_{\es}\mathbf{D}
\pmatrix{
0\cr
\delta}-\mu_u^{\T}\delta+o_p(1)\right\}
\\
&=&
\mu_{\beta}^{\T}\sum_{\es}w(\es|\widehat\delta)\mathbf{R}_{\es}\mathbf{D}
\pmatrix{ -[I,0]\mathbf{D}^{-1}\mathbf{G}_n\vspace*{2pt}\cr
-[0,I]\mathbf{D}^{-1}\mathbf{G}_n+\delta}-\mu_u^{\T}\delta+o_p(1)
\\
&=&
-\mu_{\beta}^{\T}\sum_{\es}w(\es|\widehat\delta)\mathbf{R}_{\es}\mathbf{D}
\pmatrix{
I & 0\cr
0 & 0_{\du \times \du}}\mathbf{D}^{-1}\mathbf{G}_n
\\
&&{}+
\mu_{\beta}^{\T}\sum_{\es}w(\es|\widehat\delta)
\mathbf{R}_{\es}\mathbf{D}
\pmatrix{
0\cr
\widehat\delta}-\mu_u^{\T}(\widehat\delta+[0,I]\mathbf{D}^{-1}\mathbf{G}_n)+o_p(1)
\\
&=&
-\mu_{\beta}^{\T}
\mathbf{D}^{-1}\mathbf{G}_n+\mu_{\beta}^{\T}\left\{Q(\widehat\delta)
\pmatrix{
0\cr
\widehat\delta}
-
\pmatrix{
0\cr
\widehat\delta}\right\}+o_p(1)
\\
&\stackrel{d}\longrightarrow&
-\mu_{\beta}^{\T}\mathbf{D}^{-1}\mathbf{G}+\mu_{\beta}^{\T}\left\{Q(\Delta)
\pmatrix{
0\cr
\Delta}
-
\pmatrix{
0\cr
\Delta}\right\}
\end{eqnarray*}
and thus the proof is complete.
\end{appendix}

\section*{Acknowledgments}\label{s7}
The authors would like to  thank the Co-Editors, the former Co-Editors,
one Associate Editor and three referees for their constructive comments
that substantially improved an earlier version of this paper.

\printaddresses

\end{document}